\documentclass{article}
\usepackage{geometry}
\usepackage[hidelinks]{hyperref} % hidelinks option removes the red/green boxes
\usepackage{graphicx} % Required for inserting images
\usepackage{tikz}
\usetikzlibrary{arrows.meta}
\usepackage[caption=false]{subfig} % <- Preamble
\usepackage{amsmath,amsfonts,amssymb,amsthm}
\usepackage[capitalize,nameinlink,noabbrev]{cleveref}
\newcommand{\E}{\mathcal{E}}
\newcommand{\R}{\mathbb{R}}
\newcommand{\Td}{\mathbb{T}^d}
\newcommand{\Ncal}{\mathcal{N}}

\def\dd{~\mathrm{d}}

\def\log{\mathrm{log}}
\def\sin{\mathrm{sin}}
\def\cos{\mathrm{cos}}

\def\e{\mathrm{e}}

\newtheorem{theorem}{Theorem}[section]
\newtheorem{lemma}[theorem]{Lemma}
\newtheorem{proposition}[theorem]{Proposition}
\newtheorem{remark}[theorem]{Remark}

\newtheorem{assumption}[theorem]{Assumption}
\newtheorem{example}[theorem]{Example}
\newtheorem{corollary}[theorem]{Corollary}

\numberwithin{equation}{section}

\title{Mass- and energy-preserving low-regularity integrators for nonlinear Schr\"odinger equations}
\author{Tianyu Jin\thanks{Institut f\"ur Geometrie und Praktische Mathematik, RWTH Aachen University, Im S\"usterfeld 2, 52072 Aachen, Germany\newline Email: {\texttt{jin@igpm.rwth-aachen.de}}.}
\and Georg Maierhofer\thanks{Department of Applied Mathematics and Theoretical Physics, University of Cambridge, Wilberforce Road, CB3 0WA Cambridge, United Kingdom \newline Email: {\texttt{gam37@cam.ac.uk}}.}}
\date{}

\usepackage{xcolor}

\begin{document}

\maketitle
\begin{abstract}
\noindent The nonlinear Schr\"odinger equation is a fundamental model of dispersive wave phenomena, possessing a rich geometric structure that includes the conservation of mass and energy. Preserving these conservation laws at the level of numerical discretisation remains challenging, particularly for solutions of low regularity. In particular, combining the favourable convergence properties of low-regularity integrators with the exact preservation of fundamental conservation laws has remained an open challenge. In this work, we introduce, to the best of our knowledge, the first mass- and energy-preserving low-regularity integrators for the nonlinear Schr\"odinger equation. We establish the well-posedness and convergence of the proposed schemes and demonstrate how ideas from projection methods can be used to ensure the simultaneous preservation of mass and energy. Our rigorous analysis (the first for mass- and energy-preserving schemes) is supported by extensive numerical experiments highlighting the favourable convergence and structure-preserving properties of the methods, including for rough solutions.
\end{abstract}



\section{Introduction}
In this work, we are interested in the numerical solution of the nonlinear Schr\"odinger equation (NLSE) given by
\begin{equation}
\begin{cases}
\partial_t u=i\Delta u-i\varepsilon |u|^2u,\quad t\in[0,T],\ x\in\Td,\\
u(0,x)=u_0(x),\quad x\in\Td,
\end{cases}
\label{eq:nlse}
\end{equation}
where $0<\varepsilon\leq1$ and $\mathbb T^d=(\mathbb R/2\pi\mathbb Z)^d$. The NLSE is a fundamental model of dispersive wave phenomena \cite{tao2006nonlinear}, with applications ranging from fibre optics to deep-water waves \cite{sulem2007nonlinear,tai1986observation,thomas2012nonlinear}. The analysis of rough solutions, with initial data $u_0\in H^s$ for small $s>0$, remains an active area of research. Indeed, only recently, Herr and Kwak established global well-posedness of the defocusing NLSE on $\mathbb T^2$ for arbitrary initial data in $H^s, s>0$ \cite{Herr2026}. Low-regularity phenomena in closely related NLS and Schr\"odinger-map equations also arise in the study of polygonal vortex filaments \cite{banica2009,banica2020,banica2013,banicamaierhoferschratz22,feng_maierhofer_schratz_2023} and domain-wall dynamics in ferromagnetism (cf. \cite{gamayun2019domain,misguich2019domain}).

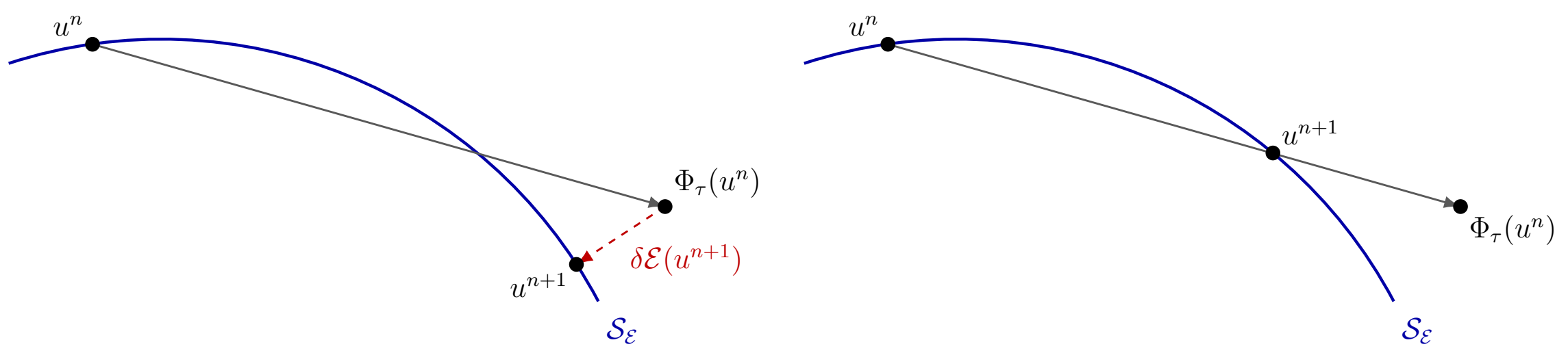
\begin{figure}[!htb]
\centering
\setlength{\abovecaptionskip}{4pt}
\subfloat[Orthogonal projection]{\resizebox{0.46\linewidth}{!}{%
\begin{tikzpicture}[
    >=Latex,
    point/.style={circle,fill=black,inner sep=2.1pt},
    level/.style={blue!65!black,very thick},
    trial/.style={gray!70!black,thick,->},
    projection/.style={red!75!black,dashed,thick,->},
    every node/.style={font=\large}
]
\begin{scope}[rotate around={-17:(1.20,-0.65)}]
    \coordinate (unA) at (-2.96,0.34);
    \coordinate (utildeA) at (5.46,0.51);
    \coordinate (uprojA) at (4.50,-0.64);
    \draw[level,shift={(0,-6)}] (125:7) arc[start angle=125,end angle=45,radius=7.0];
    \node[blue!65!black] at (5.40,-1.35) {$\mathcal{S}_{\E}$};
    \draw[trial] (unA) -- (utildeA);
    \draw[projection] (utildeA) -- node[below right,pos=0.50]
        {$\delta\E(u^{n+1})$} (uprojA);
    \node[point] at (unA) {};
    \node[above left] at (unA) {$u^n$};
    \node[point] at (utildeA) {};
    \node[above right] at (utildeA) {$\Phi_\tau(u^n)$};
    \node[point] at (uprojA) {};
    \node[below left] at (uprojA) {$u^{n+1}$};
\end{scope}
\end{tikzpicture}%
}}%
\hspace{0.02\linewidth}%
\subfloat[{Line search}]{\resizebox{0.46\linewidth}{!}{%
\begin{tikzpicture}[
    >=Latex,
    point/.style={circle,fill=black,inner sep=2.1pt},
    level/.style={blue!65!black,very thick},
    trial/.style={gray!70!black,thick,->},
    every node/.style={font=\large}
]
\begin{scope}[rotate around={-17:(1.20,-0.65)}]
    \coordinate (unB) at (-2.96,0.34);
    \coordinate (utildeB) at (5.46,0.51);
    \coordinate (ugammaB) at (2.70,0.46);
    \draw[level,shift={(0,-6)}] (125:7) arc[start angle=125,end angle=45,radius=7.0];
    \node[blue!65!black] at (5.40,-1.35) {$\mathcal{S}_{\E}$};
    \draw[trial] (unB) -- (utildeB);
    \node[point] at (unB) {};
    \node[above left] at (unB) {$u^n$};
    \node[point] at (utildeB) {};
    \node[below right] at (utildeB) {$\Phi_\tau(u^n)$};
    \node[point] at (ugammaB) {};
    \node[above right] at (ugammaB) {$u^{n+1}$};
\end{scope}
\end{tikzpicture}%
}}
\caption{Orthogonal projection and line search onto $\mathcal{S}_{\E}=\{f\,|\,\E(f) = \E(u_0)\}$.}
\label{fig:comparison_diagram}
\end{figure}

As a result, recent years have seen a significant research effort devoted to the development of low-regularity integrators, i.e. methods that can accurately and effectively approximate rough solutions under reduced regularity assumptions. This includes early work by \cite{LRI1,bai2023,ALAMABRONSARD2023114632}, as well as the development of general frameworks for low-regularity integration of dispersive nonlinear systems \cite{LRI2,rousset2021general}. At the same time, the NLSE is a Hamiltonian PDE with associated conserved energy
\begin{equation}
\E(u):=\frac12\int_{\Td}\left(|\nabla u|^2+\frac{\varepsilon}{2}|u|^4\right)\,dx,
\label{eq:energy}
\end{equation}
and conserved mass
\[
\mathcal{M}(u)=\int_{\Td} |u(t,x)|^2\dd x.
\]
This motivates the design of \textit{geometric numerical methods} that retain these conservation laws at the discrete level \cite{Hairer2006}. This has led to the development of symmetric and symplectic low-regularity integrators \cite{feng2024explicit,AlamaBronsard24, maierhofer2023bridging,alamabronsard_et_al23} that exhibit favourable energy conservation properties over long time intervals (cf.\ also \cite{faou2026}). Another interesting direction has been the study of relaxation-based approaches for the exact preservation of invariants, including recent mass-preserving low-regularity integrators \cite{li2025,akrivis2025,bai2024}, following earlier work for classical numerical schemes by \cite{ketcheson2019,Ranocha2020}. A complementary approach is provided by orthogonal projection methods \cite{Hairer2006,Celledoni2008} (cf. Figure~\ref{fig:comparison_diagram}(a)). Despite these advances, combining the favourable convergence properties of low-regularity integrators with the exact preservation of both mass and energy remains challenging. While recent relaxation-based low-regularity schemes ensure exact mass conservation, simultaneous preservation of the physical energy has remained largely unexplored in this setting, with the exception of recent work \cite{ranocha2026highordermassenergyconservingmethods} which has constructed mass- and energy-preserving schemes for Runge--Kutta base schemes on ordinary differential equations. Specifically, the work by \cite{ranocha2026highordermassenergyconservingmethods} does not include any convergence analysis for partial differential equations and, due to the CFL and regularity requirements the resulting schemes are unsuitable for the low regularity regime for the NLSE. More generally, the simultaneous preservation of different invariants is a delicate problem in geometric numerical integration. In particular, Chartier, Faou, and Murua~\cite{Chartier2006} showed that a B-series integrator satisfying the universal conditions for preserving quadratic first integrals and Hamiltonians must coincide formally with a scaled exact flow. These universal requirements do not preclude corrections tailored to selected invariants of a particular equation and it is thus a natural question to ask whether one might design methods that simultaneously preserve mass and energy.

\textbf{Key results.} In this present work, we introduce a novel low-regularity framework of relaxation-based numerical schemes for the NLSE which achieve this combined preservation of invariants. Figure~\ref{fig:comparison_diagram} illustrates two ways to enforce energy conservation. Orthogonal projection corrects along the $L^2$ energy gradient $\delta\E(u)$, whose Laplacian term requires two additional spatial derivatives to control in $H^s$. To avoid this loss, we use a line search along the baseline increment (Figure~\ref{fig:comparison_diagram}(b)). In particular, our framework combines a broad class of admissible baseline low-regularity integrators with mass normalization and scalar relaxation to preserve exactly the mass and energy of equation~\eqref{eq:nlse}, leading to the construction of a central scheme of the form
\begin{align*}
u^{n+1}
= \lambda_n \left[
u^n+\gamma_n(\Phi_\tau(u^n)-u^n)
\right],
\qquad
\lambda_n =
\sqrt{\frac{\mathcal M(u^n)}
{\mathcal M\left(u^n+\gamma_n(\Phi_\tau(u^n)-u^n)\right)}},
\end{align*}
where $\Phi_\tau$ denotes a baseline low-regularity integrator and the scalar parameter $\gamma_n$ is chosen to enforce energy conservation. The resulting scheme preserves both mass and energy exactly, with the numerical time advanced according to $t_{n+1}=t_n+\gamma_n\tau$, i.e. $u^{n+1}\approx u(t_n+\gamma_n\tau)$. For this class of methods we establish, under suitable regularity and spectral nondegeneracy assumptions, and for sufficiently small $\varepsilon$ and $\tau$, the existence and uniqueness of a relaxation parameter $\gamma_n$ close to $1$ (Theorem~\ref{thm:gamma}). Moreover, we prove that this parameter is in fact the unique nontrivial real solution of the energy conservation equation (Theorem~\ref{thm:global}). For a baseline integrator of order $p$, our analysis yields the quantitative estimates
\begin{equation}
|1-\gamma_n| \lesssim \varepsilon\tau^{p-1},
\qquad
|1-\lambda_n| \lesssim \varepsilon\tau^{p+1}.
\end{equation}
Building on these results, we establish global error estimates in Sobolev spaces (Theorem~\ref{thm:globalerror2}).

An important feature of our novel framework is its computational simplicity. The correction requires only the determination of a scalar relaxation parameter, rather than the solution of a nonlinear system. Indeed, exploiting the particular structure of the
NLSE energy functional, the conservation equation reduces, after factoring out the trivial root, to a polynomial equation of degree at most three (Proposition~\ref{prop:polynomial}). The resulting correction can thus be computed efficiently without an iterative nonlinear solver.

The remainder of this manuscript is structured as follows. In Section~\ref{sec:main_scheme} we describe the main construction of our novel methodology, together with some theoretical preliminaries. This is followed, in Section~\ref{sec:relaxation_param_theory} with the theoretical study of the solvability of the energy conservation equation with respect to the relaxation parameter and, in Section~\ref{sec:convergence_analysis}, with a detailed convergence analysis of the resulting schemes in the low-regularity regime. Finally, we evaluate our novel method and benchmark against prior work in Section~\ref{sec:numerics}, before providing concluding remarks and an outlook for future work in Section~\ref{sec:conclusions}.

\section{Construction of mass- and energy-preserving integrators}\label{sec:main_scheme}
In the following, we will describe preliminary results and the construction of our novel low-regularity integrators, before studying their properties rigorously in Sections~\ref{sec:relaxation_param_theory} and \ref{sec:convergence_analysis}.
\subsection{Preliminaries and notations}
Throughout, we will follow the convention that integrals and $L^q$ norms on $\mathbb T^d=(\mathbb R/2\pi\mathbb Z)^d$
use the normalized Haar measure, so $\int_{\mathbb T^d}1\,dx=1$. We then use the convention that the Fourier coefficients for $f\in L^2(\mathbb{T}^d)$ are given by
$$
\hat{f}_{k  
}=\int_{\mathbb T^d}f(x)\e^{-ik\cdot x}\,dx,
$$
such that, in particular,
Parseval's identity reads $\|f\|_{L^2}^2=\sum_k|\widehat f_k|^2$. Using this convention we can define Sobolev norms
\begin{align*}
    \|f\|_{H^s}^2:=\sum_{k\in\mathbb{Z}^d}\langle k\rangle^{2s}|\hat{f}_k|^2, \quad \text{where } \langle k\rangle:=\left\{ \begin{array}{cc}
     |k|\quad&\text{if } k\neq0,\\
      1\quad&\text{if } k=0.\\
\end{array}
\right.
\end{align*}
In subsequent chapters, we will analyse our methods in the periodic Sobolev spaces given by
\begin{align*}
    H^{s}(\Omega):=\{f\in L^2(\Omega)\,\vert\, \|f\|_{H^s}<\infty\}.
\end{align*}
Moreover, in the interest of notational simplicity, for $1\leq q \leq\infty$, we write
\[
\|\cdot\|_q := \|\cdot\|_{L^q(\Td)},
\]
and we employ the notation $A\lesssim B$ means that $A\leq CB$ to mean that (unless stated otherwise) $C>0$ denotes a constant independent of $\varepsilon$ and $\tau$, whose value may change between occurrences. Constants may depend on the solution norms indicated in each result. Independence of $\varepsilon$ and $\tau$ is understood uniformly whenever these norms are uniformly bounded.

In the analysis of later sections, the local Lipschitz continuity of the energy functional \eqref{eq:energy} will play a crucial role. In particular we will use the following result.
\begin{lemma}\label[lemma]{lem:energy}
    Let $s>d/2$ and $s\geq 1$. For any $v,w\in H^s$ and $0<\varepsilon\leq1$, there is a constant $C_E>0$ that depends on $s,d,\|v\|_{H^s}$ and $\|w\|_{H^s}$, and is independent of $\varepsilon$, such that
    \[
    |\E(v)-\E(w)|\leq C_E\|v-w\|_{H^s}.
    \]
\end{lemma}
\begin{proof}
From \eqref{eq:energy} we know that
    \begin{align*}
        |\E(v)-\E(w)|
        &= \left|\frac12(\|\nabla v\|_{L^2}^2-\|\nabla w\|_{L^2}^2)
        +\frac{\varepsilon}{4}\int_{\mathbb T^d} (|v|^4-|w|^4)\dd x\right|\\
        &\leq \frac12\left|\operatorname{Re}\langle\nabla(v-w),\nabla(v+w)\rangle\right|
        +\frac{\varepsilon}{4}\int_{\mathbb T^d}\left||v|^4-|w|^4\right|\dd x.
    \end{align*}
    We then estimate the two terms respectively. Obviously,
    \begin{align*}
        \frac12\left|\operatorname{Re}\langle\nabla(v-w),\nabla(v+w)\rangle\right|
        \leq \frac12(\|v\|_{H^s}+\|w\|_{H^s})\|v-w\|_{H^s}.
    \end{align*}
    Moreover,
    \begin{align*}
        \left||v|^4-|w|^4\right|
        &=(|v|^3+|v|^2|w|+|v||w|^2+|w|^3)\left||v|-|w|\right|\\
        &\leq 2(|v|^3+|w|^3)|v-w|.
    \end{align*}
    Since $s>d/2$, the Sobolev embedding $H^s(\mathbb T^d)\hookrightarrow L^\infty(\mathbb T^d)$ gives
    \begin{align*}
        \int_{\mathbb T^d}\left||v|^4-|w|^4\right|\dd x
        &\leq 2(\|v\|_{L^\infty}^3+\|w\|_{L^\infty}^3)\|v-w\|_{L^1}\\
        &\leq C_0(\|v\|_{H^s}^3+\|w\|_{H^s}^3)\|v-w\|_{H^s}.
    \end{align*}
    Therefore, using $0<\varepsilon\leq1$, there is a constant $C_E>0$ depending only on $s,d,\|v\|_{H^s}$ and $\|w\|_{H^s}$ such that
    \[
    |\E(v)-\E(w)|\leq C_E\|v-w\|_{H^s}.
    \]
\end{proof}
Finally, throughout the following sections, we denote by $\phi_t(v)$ the exact flow of \eqref{eq:nlse} with initial value $v$.

\subsection{Key assumptions} The basis of our construction of mass- and energy-preserving schemes are low-regularity integrators that satisfy the following assumptions.
\begin{assumption}\label[assumption]{assump:lri}
For any $v\in H^s$, $s>d/2$, the low-regularity integrators for \eqref{eq:nlse} with time stepsize $\tau$, $\Phi_{\tau}:H^s\rightarrow H^s$, can be written as
\[
\Phi_\tau(v)=e^{i\tau\Delta}v+\Ncal_\tau(v),
\]
and they have the following properties.

\emph{(i)} For any $\tau>0$, $\Ncal_\tau(0)=0$.

\emph{(ii)} There is a continuous function $M:\R_{\ge0}\times\R_{\ge0}\to\R_{\ge0}$ such that for any $v,w\in H^s$ and $0<\varepsilon,\tau\le1$, we have
\[
\|\Ncal_\tau(v)-\Ncal_\tau(w)\|_{H^s}
\le \varepsilon\tau M(\|v\|_{H^s},\|w\|_{H^s})\|v-w\|_{H^s}.
\]

\emph{(iii)} There are constants $p\ge1$, $\eta\ge0$ and a continuous, increasing function $C:\R_{\ge0}\to\R_{\ge0}$ such that the local truncation error satisfies
\[
\|\Phi_\tau(u(t))-u(t+\tau)\|_{H^s}
=\|\Ncal_\tau(u(t))-I_\tau(u(t))\|_{H^s}
\le \varepsilon\tau^{p+1}C\left(\sup_{\xi\in[t,t+\tau]}\|u(\xi)\|_{H^{s+\eta}}\right),
\]
for any $0<\varepsilon,\tau\le1$ and $t\ge0$ with $t+\tau\le T$. Here, $I_\tau(v):=-i\varepsilon e^{i\tau\Delta}\int_0^\tau e^{-is\Delta}\bigl(|\phi_s(v)|^2\phi_s(v)\bigr)\,ds.$
\end{assumption}
These assumptions are, in fact, sufficiently general to cover the majority of existing low-regularity integrators for the nonlinear Schr\"odinger equation. In particular, in the following examples we present two widely-used low-regularity schemes for the NLSE \eqref{eq:nlse} that satisfy the above assumption.
\begin{example}\label[example]{eg:nlse_lri}
    (i) \textit{First-order resonance-based low-regularity integrator (LRI1).} This was first proposed in \cite{LRI1} with the following formula
    \begin{align}\label{eq:nlse_lri1}
        \Phi_\tau^{LRI1}(v) = \e^{i\tau\Delta}v-i\varepsilon\tau \e^{i\tau\Delta} \left[v^2\varphi_{1}(-2i\tau\Delta) \overline{v}\right],
    \end{align}
    where $\varphi_1(z) = \frac{\e^z-1}{z}$ for $z\in\mathbb{C}$. In this example, $s>\frac{d}{2}$ and $\eta=1$.

    \noindent
    (ii) \textit{Second-order resonance-based low-regularity integrator (LRI2).} This was introduced by Bruned \& Schratz \cite{LRI2} with the following formula
    \begin{align}\label{eq:nlse_lri2}
        \Phi_\tau^{LRI2}(v) &= \e^{i\tau\Delta}v-i\varepsilon\tau \e^{i\tau\Delta} (v^2[\varphi_{1}(-2i\tau\Delta)-\varphi_2(-2i\tau\Delta)] \overline{v})\\\nonumber
        &\quad -i\varepsilon\tau(\e^{i\tau\Delta}v)^2[\e^{i\tau\Delta}\varphi_2(-2i\tau\Delta)\overline{v}] - \varepsilon^2\frac{\tau^2}{2}\e^{i\tau\Delta}(|v|^4v),
    \end{align}
    where $\varphi_2(z) = \frac{\e^z-\varphi_1(z)}{z}$ for $z\in\mathbb{C}$. In this example, $s>\frac{d}{2}$ and $\eta=2$.
\end{example}
Based on the above assumptions we can start by considering the mass and energy defects after applying $\Phi_\tau$ to an initial value $v$,
\[
\delta_M(v):=\|\Phi_\tau(v)\|_2^2-\|v\|_2^2,
\qquad
\delta_E(v):=\E(\Phi_\tau(v))-\E(v).
\]
The properties in Assumption~\ref{assump:lri} allow us to derive estimates for these defects. 
\begin{corollary}\label[corollary]{lem:defects}
Let $s>\max\{1,d/2\}$, $0<\varepsilon,\tau\leq1$ and $v\in H^{s+\eta}$. Suppose that $\Phi_\tau$ satisfies Assumption~\ref{assump:lri} and that $\phi_t(v)\in H^{s+\eta}$ for $0\leq t\leq\tau$. Then
\begin{equation}
|\delta_M(v)|\leq C\varepsilon\tau^{p+1}, \quad|\delta_E(v)|\leq C\varepsilon\tau^{p+1}.
\label{eq:defectbounds2}
\end{equation}
The constant $C>0$ depends on $s,d$ and $\sup_{0\leq t\leq\tau}\|\phi_t(v)\|_{H^{s+\eta}}$.
\end{corollary}
\begin{proof}
From Assumption~\ref{assump:lri} (iii) we know that
\[
\|\Phi_\tau(v)-\phi_\tau(v)\|_{H^s}
\lesssim\varepsilon\tau^{p+1},
\]
with an implicit constant that depends on $s,d$ and $\sup_{0\leq t\leq\tau}\|\phi_t(v)\|_{H^{s+\eta}}$. Using the exact conservation laws and Lemma~\ref{lem:energy},
\begin{align*}
|\delta_M(v)|
&\le\bigl(\|\Phi_\tau(v)\|_2+\|\phi_\tau(v)\|_2\bigr)
\|\Phi_\tau(v)-\phi_\tau(v)\|_2
\lesssim\varepsilon\tau^{p+1},\\
|\delta_E(v)|
&=|\E(\Phi_\tau(v))-\E(\phi_\tau(v))|
\lesssim\|\Phi_\tau(v)-\phi_\tau(v)\|_{H^s}
\lesssim\varepsilon\tau^{p+1}.
\end{align*}
\end{proof}
For the purpose of our analysis, we will need to find an upper bound for the increment $\|\Phi_\tau(v)-v\|_{H^s}$, which is possible with additional regularity requirements.
\begin{lemma}\label[lemma]{lem:incrementupper}
    Assume that $v\in H^{s+2}$, $s>\frac{d}{2}$, and that $\Phi_\tau$ satisfies Assumption~\ref{assump:lri}. Then we have that
    \[
    \|\Phi_\tau(v)-v\|_{H^s}\leq C\tau,\qquad 0<\varepsilon,\tau\leq1,
    \]
    where $C>0$ depends on $s,d$, and $\|v\|_{H^{s+2}}$.
\end{lemma}
\begin{proof}
    Using the formula for the baseline integrator, we know that
    \[
    \|\Phi_\tau(v)-v\|_{H^s}\leq \|(e^{i\tau\Delta}-I)v\|_{H^s} + \|\mathcal{N}_\tau(v)\|_{H^s}.
    \]
    By Parseval's identity and $|e^{ix}-1|\le |x|$,
    \begin{align*}
    \|(e^{i\tau\Delta}-I)v\|_{H^s}^2
    &=\sum_{k\in\mathbb Z^d}
    \langle k\rangle^{2s}
    |e^{-i\tau|k|^2}-1|^2|\widehat v_k|^2\\
    &\le \tau^2\sum_{k\in\mathbb Z^d}
    \langle k\rangle^{2s}|k|^4|\widehat v_k|^2\\
    &\le \tau^2\|v\|_{H^{s+2}}^2.
    \end{align*}
    Consequently,
    \[
    \|(e^{i\tau\Delta}-I)v\|_{H^s}
    \le \tau\|v\|_{H^{s+2}}.
    \]
    From Assumption~\ref{assump:lri} (i) and (ii), we know that 
    \begin{equation*}
        \|\mathcal{N}_\tau(v)\|_{H^s}\lesssim \varepsilon\tau.
    \end{equation*}
    Therefore, we obtain that
    \[
    \|\Phi_\tau(v)-v\|_{H^s}\lesssim \tau.
    \]
\end{proof}

\subsection{Novel mass and energy conserving low-regularity integrators} 
Throughout the rest of this paper, we assume that $\Phi_\tau$ satisfies Assumption~\ref{assump:lri}. Starting from the previous approximation $u^n\approx u(t_n)$, we propose the following time-relaxation scheme based on $\Phi_\tau$
\begin{equation}
\begin{aligned}
u^{n+1}&=\lambda_n\bigl[u^n+\gamma_n (\Phi_\tau(u^n)-u^n)\bigr],\qquad
\lambda_n=\sqrt{\frac{\mathcal M\left(u^n\right)}
{\mathcal M\left(u^n+\gamma_n (\Phi_\tau(u^n)-u^n)\right)}}.
\end{aligned}
\label{eq:scheme2}
\end{equation}
$u^{n+1}\approx u(t_{n+1}):=u(t_n + \gamma_n\tau)$.
At the current time step, the relaxation parameter $\gamma_n$ is chosen to enforce energy conservation. Define
\begin{equation}
H(z):=\E\left(
\frac{\|u^n\|_{2}\,[u^n+z(\Phi_\tau(u^n)-u^n)]}
{\|u^n+z(\Phi_\tau(u^n)-u^n)\|_{2}}
\right)-\E(u^n),
\label{eq:H}
\end{equation}
we will show in the following that it is possible to find a root $\gamma_n$ close to $1$ such that
\[
H(\gamma_n) = 0.
\]
By such a construction, we then immediately have
\[
\mathcal M\left(u^{n+1}\right) = \mathcal M\left(u^n\right) \quad\text{and}\quad \E(u^{n+1}) = \E(u^n).
\]
The update \eqref{eq:scheme2} builds on time relaxation methods \cite{ketcheson2019,Ranocha2020}, in which a scalar parameter rescales both the numerical increment and the time step. Our scheme simplifies the quadratic-preserving relaxation method of \cite{ranocha2026highordermassenergyconservingmethods}. Their method first normalizes $\Phi_\tau(u^n)$ to have the same mass as $u^n$ before forming the relaxation increment. We form the increment directly from $\Phi_\tau(u^n)$ and normalize only the relaxed trial state. The two updates coincide when the baseline integrator is already mass-preserving.

For $z\in\R$, we define $D(z)$ as the mass of the unnormalized trial state with the baseline increment scaled by $z$, i.e.
\[
D(z):=\mathcal M\left(u^n+z(\Phi_\tau(u^n)-u^n)\right)=\|u^n+z(\Phi_\tau(u^n)-u^n)\|_2^2.
\]
To conclude this section, we derive some algebraic identities that is related to the energy residual $H$, which will play a role in later analysis.
\begin{lemma}[Algebraic identities]\label[lemma]{lem:identities2}
For every $z\in\mathbb R$,
\begin{align}
D(z)&=\|u^n\|_2^2+z(z-1)\|\Phi_\tau(u^n)-u^n\|_2^2+z\delta_M,
\label{eq:massidentity}\\
\|\nabla(u^n+z(\Phi_\tau(u^n)-u^n))\|_2^2
&=\|\nabla u^n\|_2^2+z(z-1)\|\nabla (\Phi_\tau(u^n)-u^n)\|_2^2\nonumber\\
&\quad+z\bigl(\|\nabla\Phi_\tau(u^n)\|_2^2-\|\nabla u^n\|_2^2\bigr).
\label{eq:kineticidentity2}
\end{align}
Where $D(z)>0$, the energy residual is
\begin{equation}
H(z)=\frac{\|u^n\|_2^2\|\nabla(u^n+z(\Phi_\tau(u^n)-u^n))\|_2^2}{2D(z)}
+\frac{\varepsilon\|u^n\|_2^4\|u^n+z(\Phi_\tau(u^n)-u^n)\|_4^4}{4D(z)^2}-\E(u^n).
\label{eq:Hexpanded1}
\end{equation}
\end{lemma}
\begin{proof}
For any two vectors $v,w$ in a complex Hilbert space,
\begin{align*}
\|v+zw\|^2
&=\|v\|^2+2z\operatorname{Re}\langle v,w\rangle+z^2\|w\|^2\\
&=\|v\|^2+z\bigl(\|v+w\|^2-\|v\|^2-\|w\|^2\bigr)+z^2\|w\|^2\\
&=\|v\|^2+z(z-1)\|w\|^2+z\bigl(\|v+w\|^2-\|v\|^2\bigr).
\end{align*}
Plugging in $v = u^n, w = \Phi_\tau(u^n)-u^n$ gives \eqref{eq:massidentity}, while plugging in $v = \nabla u^n,w = \nabla (\Phi_\tau(u^n)-u^n)$ gives \eqref{eq:kineticidentity2}.
Substitution of the normalization factor into \eqref{eq:energy} gives
\eqref{eq:Hexpanded1}.
\end{proof}

\section{Existence and uniqueness of the relaxation parameter}\label{sec:relaxation_param_theory}
For our proposed scheme to work, we first need to know whether the algebraic equation $H(z)=0$ has a solution close to $1$. In this section, we study the existence and uniqueness of the nonzero root of $H$ for sufficiently small $\tau$ and $\varepsilon$. To begin with, we provide some estimates for $H$ and its derivatives. Throughout the remainder of this work, we assume that $u^n$ is nonzero, that is, $\|u^n\|_2>0$.
\subsection{Estimates for the energy residual}
Exploiting the identity \eqref{eq:Hexpanded1}, we obtain the following estimate of the absolute value of $H(1)$.
\begin{corollary}\label[corollary]{lem:H1}
Assume $u^n\in H^{s+\max\{2,\eta\}}$, $s>\max\{1,d/2\}$, and $\phi_t(u^n)\in H^{s+\max\{2,\eta\}}$ for $0\leq t\leq\tau$. Then $H(0)=0$. For sufficiently small $\tau$,
\begin{equation}
|H(1)|\leq C\varepsilon\tau^{p+1}.
\label{eq:H1}
\end{equation}
The constant $C>0$ %and the step-size threshold 
depends on $s,d,\|u^n\|_2$, and $\sup_{0\leq t\leq\tau}\|\phi_t(u^n)\|_{H^{s+\max\{2,\eta\}}}$.
\end{corollary}
\begin{proof}
Independent of $\varepsilon$, we know for a sufficiently small $\tau$,
\[
|\delta_M|\le\tfrac12\|u^n\|_2^2,
\qquad D(1)=\|u^n\|_2^2+\delta_M\ge\tfrac12\|u^n\|_2^2.
\]
The energy formula gives
\begin{align*}
H(1)
&=\delta_E-\frac{\delta_M\|\nabla\Phi_\tau(u^n)\|_2^2}{2D(1)}-\frac{\varepsilon\delta_M(\|u^n\|_2^2+D(1))
\|\Phi_\tau(u^n)\|_4^4}{4D(1)^2}.
\end{align*}
Thus, Corollary~\ref{lem:defects} implies that
\[
|H(1)|\lesssim |\delta_E|+|\delta_M|
\lesssim\varepsilon\tau^{p+1}.
\]
\end{proof}

To establish the existence of a root of $H$ close to $1$, we complement the estimate for $H(1)$ with a lower bound for $H'(1)$. To identify its leading contribution for small $\varepsilon$, we fix $u^n$ and set $\varepsilon=0$ in both the baseline step and the energy. The resulting energy residual is denoted by
\[
H_{l}(z)=\frac12\left(\frac{\|u^n\|_2^2\,\|\nabla(u^n+z (e^{i\tau\Delta}-I)u^n)\|_2^2}{\|u^n+z (e^{i\tau\Delta}-I)u^n\|_2^2}-\|\nabla u^n\|_2^2\right).
\]
Define the quantity $\Sigma(\tau)$ as the slope of the linear energy residual $H_l$ at $1$, i.e.
\begin{align}
\Sigma(\tau)
&:=H_{l}'(1)\nonumber\\
&=\frac12\left(
\|\nabla (e^{i\tau\Delta}-I)u^n\|_2^2
-\frac{\|\nabla u^n\|_2^2}{\|u^n\|_2^2}\|(e^{i\tau\Delta}-I)u^n\|_2^2
\right)\nonumber\\
&=\sum_k \bigl(1-\cos(\tau|k|^2)\bigr)
\left(|k|^2-\frac{\|\nabla u^n\|_2^2}{\|u^n\|_2^2}\right)
|\widehat u_k^n|^2.\label{eq:Sigma}
\end{align}
The last line comes from Parseval’s identity. 
%We assume that the Fourier support contains two different radii. More precisely, fix $N\in\mathbb N$ such that
%\begin{equation}
%\sigma_N:=\frac{1}{\pi^2\|u^n\|_2^2}
%\sum_{|k|,|l|\le N}(|k|^2+|l|^2)(|k|^2-|l|^2)^2
%|\widehat u^n_k|^2|\widehat u^n_l|^2>0.
%\label{eq:sigmaN}
%\end{equation}
We next establish a positive lower bound for $\Sigma(\tau)$ under a spectral nondegeneracy condition.
\begin{proposition}\label[proposition]{lem:Sigmalower}
Assume $u^n\in H^s$, $s>\frac32$. For $N\in\mathbb N$, define
\begin{equation}\label{eq:sigmaN}
\sigma_N:=\frac{1}{\pi^2\|u^n\|_2^2}
\sum_{|k|,|l|\le N}(|k|^2+|l|^2)(|k|^2-|l|^2)^2
|\widehat u^n_k|^2|\widehat u^n_l|^2>0.
\end{equation}
Then for any fixed $N$ satisfying $\sigma_N>0$, there is a threshold $0<\tau_1\le\pi/(2N^2)$ such that, for $0<\tau\le\tau_1$,
\begin{equation}
\Sigma(\tau)\ge\tfrac12\sigma_N\tau^2.
\label{eq:SigmamuLower}
\end{equation}
The threshold $\tau_1$ depends on $s,d,N,\sigma_N,\|u^n\|_{H^s}$, and $\|u^n\|_2$.
\end{proposition}

\begin{proof}
\emph{Step 1:} We first study the low-frequency contribution to $\Sigma(\tau)$.
The definition of $\Sigma(\tau)$ \eqref{eq:Sigma} gives
\begin{align*}
    \Sigma(\tau)
    &=\left(\sum_k\bigl(1-\cos(\tau|k|^2)\bigr)|k|^2|\hat{u}_k^n|^2\right)-\frac{\sum_l|l|^2|\hat{u}_l^n|^2}{\|u^n\|_2^2}\sum_k\bigl(1-\cos(\tau|k|^2)\bigr)|\hat{u}_k^n|^2\\
    &=\frac1{\|u^n\|_2^2}\sum_{k,l}
    \bigl(1-\cos(\tau|k|^2)\bigr)
    (|k|^2-|l|^2)
    |\widehat u_k^n|^2|\widehat u_l^n|^2.
\end{align*}
Interchanging $k$ and $l$ gives
\[
\Sigma(\tau)=-\frac1{\|u^n\|_2^2}\sum_{k,l}
    \bigl(1-\cos(\tau|l|^2)\bigr)
    (|k|^2-|l|^2)
    |\widehat u_k^n|^2|\widehat u_l^n|^2.
\]
Hence, we have that
\begin{equation}
    \Sigma(\tau)=\frac1{2\|u^n\|_2^2}\sum_{k,l}
    \bigl(\cos(\tau|l|^2)-\cos(\tau|k|^2)\bigr)
    (|k|^2-|l|^2)
    |\widehat u_k^n|^2|\widehat u_l^n|^2.
    \label{eq:pairform}
\end{equation}
We know that $\Sigma(\tau)$ is absolutely convergent, since
\[
\Sigma(\tau)\leq\frac1{2\|u^n\|_2^2}
\sum_{k,l}2(|k|^2+|l|^2)
|\widehat u_k^n|^2|\widehat u_l^n|^2
=2\|\nabla u^n\|_2^2<\infty.
\]

Define two sets $A_\tau=\{k:\tau|k|^2\le\pi\}, \,B_\tau=\mathbb Z^d\setminus A_\tau$.
We know that for all $k,l\in A_\tau$, $\bigl(\cos(\tau|l|^2)-\cos(\tau|k|^2)\bigr)(|k|^2-|l|^2)\ge0$.
Moreover, from the fact that
\[
\cos b-\cos a
=2\sin\frac{a+b}{2}\sin\frac{a-b}{2}
\ge\frac2{\pi^2}(a+b)(a-b), \quad 0\le b\le a\le\pi/2,
\]
the contribution of $|k|,|l|\le N$ to \eqref{eq:pairform} is at least
\begin{align*}
&\frac{\tau^2}{\pi^2\|u^n\|_2^2}
\sum_{|k|,|l|\le N}
(|k|^2+|l|^2)(|k|^2-|l|^2)^2
|\widehat u_k^n|^2|\widehat u_l^n|^2=\sigma_N\tau^2,
\end{align*}
since $\tau N^2\le\pi/2$.

\emph{Step 2:} Next, we control the high-frequency contributions.
Set $R=\|u^n\|_{H^s}$. Using $s>\frac32$, we obtain
\begin{align*}
\sum_{l\in A_\tau}
(1-\cos(\tau|l|^2))|\widehat u_l^n|^2
&\lesssim
\tau^{\min\{s,2\}}
\sum_l |l|^{2\min\{s,2\}}|\widehat u_l^n|^2
\lesssim R^2\tau^{\min\{s,2\}},\\
\sum_{k\in B_\tau}|k|^2|\widehat u_k^n|^2
&=\sum_{k\in B_\tau}
|k|^{2s}|k|^{-2(s-1)}|\widehat u_k^n|^2
\le \frac{R^2}{\pi^{s-1}}\tau^{s-1},\\
\sum_{k\in B_\tau}|\widehat u_k^n|^2
&=\sum_{k\in B_\tau}
|k|^{2s}|k|^{-2s}|\widehat u_k^n|^2
\le \frac{R^2}{\pi^s}\tau^s.
\end{align*}

For $k\in B_\tau$ and $l\in A_\tau$,
\[
\cos(\tau|l|^2)-\cos(\tau|k|^2)
\ge -(1-\cos(\tau|l|^2)).
\]
Since $|k|^2>|l|^2$, it follows that
\[
\bigl(\cos(\tau|l|^2)-\cos(\tau|k|^2)\bigr)
(|k|^2-|l|^2)
\ge -|k|^2(1-\cos(\tau|l|^2)).
\]
Thus the combined contribution of
$B_\tau\times A_\tau$ and $A_\tau\times B_\tau$
is bounded below by
\begin{align*}
&-\frac1{\|u^n\|_2^2}
\left(\sum_{k\in B_\tau}|k|^2|\widehat u_k^n|^2\right)
\left(\sum_{l\in A_\tau}
(1-\cos(\tau|l|^2))|\widehat u_l^n|^2\right)\ge -C\tau^{s-1+\min\{s,2\}}.
\end{align*}
The absolute contribution of $B_\tau\times B_\tau$
is at most
\begin{align*}
&\frac2{\|u^n\|_2^2}
\left(\sum_{k\in B_\tau}|k|^2|\widehat u_k^n|^2\right)
\left(\sum_{l\in B_\tau}|\widehat u_l^n|^2\right)
\le C\tau^{2s-1},
\end{align*}
where $C$ depends on $s,R,\|u^n\|_2$.

Combining these estimates with Step 1 gives
\begin{align*}
\Sigma(\tau)
&\ge \sigma_N\tau^2
-C\left(\tau^{s-1+\min\{s,2\}}+\tau^{2s-1}\right)\\
&=\tau^2\left[
\sigma_N-C\left(
\tau^{s-3+\min\{s,2\}}+\tau^{2s-3}
\right)\right].
\end{align*}
Both exponents in the last bracket are positive because $s>3/2$.
Therefore, for sufficiently small $\tau>0$,
\[
\Sigma(\tau)\ge \tfrac12\sigma_N\tau^2.
\]
\end{proof}

\begin{remark}
    The spectral nondegeneracy condition $\sigma_N>0$ holds when there exist $k,l$ with $|k|, |l|\leq N, |k|\neq |l|$, and $\hat{u}_k^n, \hat{u}_l^n\neq0$. In practice, it suffices to choose $N$ large enough to include two nonzero Fourier coefficients whose wave numbers have different lengths.

    For $N\in\mathbb N$ and $0\ne v\in L^2(\mathbb T^d)$, define
    \begin{equation*}
    \sigma_N(v):=\frac{1}{\pi^2\|v\|_2^2}
    \sum_{|k|,|l|\le N}
    (|k|^2+|l|^2)(|k|^2-|l|^2)^2
    |\widehat v_k|^2|\widehat v_l|^2.
    \end{equation*}
    For all the one-step analysis starting from the datum $u^n$, we fix $N$ such that $\sigma_N(u^n)>0$ and write $\sigma_N=\sigma_N(u^n)$.
\end{remark}

We next estimate $H'(1)$ with the aid of $\Sigma(\tau)$ and bound $H''$ to control the behaviour of $H$ near $1$.
\begin{lemma}\label[lemma]{lem:Hprime}
Assume $u^n\in H^{s+\max\{2,\eta\}}$, $s>\max\{1,d/2\}$, and $\phi_t(u^n)\in H^{s+\max\{2,\eta\}}$ for $0\leq t\leq\tau$. For all sufficiently
small $\tau>0$, $H$ is defined on $[-2,2]$ and
\begin{equation}
|H'(1)-\Sigma(\tau)|\leq C\varepsilon\tau^2,
\qquad
\sup_{|z|\le2}|H''(z)|\leq C\tau^2.
\label{eq:Hderiv}
\end{equation}
The constant $C>0$ and the smallness threshold of $\tau$ depend on $s,d,\|u^n\|_2$, and
$\sup_{0\leq t\leq\tau}\|\phi_t(u^n)\|_{H^{s+\max\{2,\eta\}}}$.
\end{lemma}
\begin{proof}
Recall that
\begin{align*}
D(z)&=\|u^n+z(\Phi_\tau(u^n)-u^n)\|_2^2.\\
\end{align*}
Define
\begin{align*}
G(z)=\|\nabla[u^n+z(\Phi_\tau(u^n)-u^n)]\|_2^2,\quad\text{and}\quad F(z)=\|u^n+z(\Phi_\tau(u^n)-u^n)\|_4^4,
\end{align*}
then from \eqref{eq:Hexpanded1} we know that
\begin{equation}\label{eq:Hexpanded2}
    H(z)=\frac{\|u^n\|_2^2G(z)}{2D(z)}
+\frac{\varepsilon\|u^n\|_2^4F(z)}{4D(z)^2}-E(u^n).
\end{equation}
Differentiating \eqref{eq:Hexpanded2} once, we obtain that
\begin{equation*}
H'(z)=\frac{\|u^n\|_2^2G'(z)}{2D(z)}
+\frac{\varepsilon\|u^n\|_2^4F'(z)}{4D(z)^2}-
\left[(2z-1)\|\Phi_\tau(u^n)-u^n\|_2^2+\delta_M\right]
\left(
\frac{\|u^n\|_2^2G(z)}{2D(z)^2}
+\frac{\varepsilon\|u^n\|_2^4F(z)}{2D(z)^3}
\right).
\end{equation*}
In particular, at $z=1$,
\begin{equation}\label{eq:Hderiv1}
    H'(1)=\frac{\|u^n\|_2^2G'(1)}{2D(1)}
+\frac{\varepsilon\|u^n\|_2^4F'(1)}{4D(1)^2}-
\bigl(\|\Phi_\tau(u^n)-u^n\|_2^2+\delta_M\bigr)
\left(
\frac{\|u^n\|_2^2G(1)}{2D(1)^2}
+\frac{\varepsilon\|u^n\|_2^4F(1)}{2D(1)^3}
\right),
\end{equation}
where $D(1)=\|u^n\|_2^2+\delta_M$.

Differentiating \eqref{eq:Hexpanded2} twice gives
\begin{align}\label{eq:Hderiv2}
H''(z)=\frac{\|u^n\|_2^2}{2}\left(
\frac{G''}{D}
-\frac{2G'D'+GD''}{D^2}
+\frac{2G(D')^2}{D^3}
\right)+\frac{\varepsilon\|u^n\|_2^4}{4}\left(
\frac{F''}{D^2}
-\frac{4F'D'+2FD''}{D^3}
+\frac{6F(D')^2}{D^4}
\right).
\end{align}

Lemma~\ref{lem:incrementupper} and Corollary~\ref{lem:defects} give
\[
\|\Phi_\tau(u^n)-u^n\|_{H^s}\lesssim\tau,
\qquad
|\delta_M|+|\delta_E|
\lesssim\varepsilon\tau^{p+1}
\lesssim\varepsilon\tau^2,
\]
where $p\ge1$ and $\tau\le1$.

\emph{Step 1:} We first show the positivity of $D(z)$ and a bound on $H''(z)$ when $|z|\leq2$.
The identity
\[
D(z)=\|u^n\|_2^2+z\delta_M
+z(z-1)\|\Phi_\tau(u^n)-u^n\|_2^2
\]
implies that for $|z|\le2$,
\[
D(z)\ge
\|u^n\|_2^2-2|\delta_M|
-6\|\Phi_\tau(u^n)-u^n\|_2^2
\ge\frac12\|u^n\|_2^2
\]
for sufficiently small $\tau$. Hence $H$ is defined on $[-2,2$].
Moreover,
\[
D'(z)=\delta_M+(2z-1)\|\Phi_\tau(u^n)-u^n\|_2^2,
\qquad
D''(z)=2\|\Phi_\tau(u^n)-u^n\|_2^2,
\]
and therefore
\[
|D'(z)|+|D''(z)|\lesssim\tau^2.
\]
Since
\[
\|\Phi_\tau(u^n)-u^n\|_2
\le\|\Phi_\tau(u^n)-u^n\|_{H^s}
\lesssim\tau,
\]
and
\[
\|\Phi_\tau(u^n)-u^n\|_{H^s}
\le2\|u^n\|_{H^s}+\|\mathcal N_\tau(u^n)\|_{H^s}
\lesssim1,
\]
using $H^s\hookrightarrow L^\infty$, we know that for $|z|\leq2$,
\begin{align*}
\|u^n+z(\Phi_\tau(u^n)-u^n)\|_\infty
&\lesssim
\|u^n\|_{H^s}
+2\|\Phi_\tau(u^n)-u^n\|_{H^s}
\lesssim1,\\
\|u^n+z(\Phi_\tau(u^n)-u^n)\|_{H^s}
&\le
\|u^n\|_{H^s}
+2\|\Phi_\tau(u^n)-u^n\|_{H^s}
\lesssim1.
\end{align*}

For $G$ and its derivatives, the Cauchy--Schwarz inequality yields
\begin{align*}
|G(z)|
&=\|\nabla[u^n+z(\Phi_\tau(u^n)-u^n)]\|_2^2
\lesssim1,\\
|G'(z)|
&\le
2\|\nabla[u^n+z(\Phi_\tau(u^n)-u^n)]\|_2
\|\nabla(\Phi_\tau(u^n)-u^n)\|_2
\lesssim\tau,\\
|G''(z)|
&=2\|\nabla(\Phi_\tau(u^n)-u^n)\|_2^2
\lesssim\tau^2.
\end{align*}
For $F$ and its derivatives, H\"older's inequality gives
\begin{align*}
|F(z)|
&\le
\|u^n+z(\Phi_\tau(u^n)-u^n)\|_\infty^2
\|u^n+z(\Phi_\tau(u^n)-u^n)\|_2^2
\lesssim1,\\
|F'(z)|
&\le
4\|u^n+z(\Phi_\tau(u^n)-u^n)\|_\infty^2
\|u^n+z(\Phi_\tau(u^n)-u^n)\|_2
\|\Phi_\tau(u^n)-u^n\|_2
\lesssim\tau,\\
|F''(z)|
&\le
12\|u^n+z(\Phi_\tau(u^n)-u^n)\|_\infty^2
\|\Phi_\tau(u^n)-u^n\|_2^2
\lesssim\tau^2.
\end{align*}
Therefore, uniformly for $|z|\le2$,
\[
|G(z)|+|F(z)|\lesssim1,\qquad
|G'(z)|+|F'(z)|\lesssim\tau,\qquad
|G''(z)|+|F''(z)|\lesssim\tau^2.
\]
All functions on the right hand side of \eqref{eq:Hderiv2} are evaluated at $z$. Consequently,
\[
\sup_{|z|\le2}|H''(z)|\lesssim\tau^2.
\]

\emph{Step 2:} We then turn to the estimate of $H'(1)$.
Define the unnormalized energy residual
\[
P(z)=E(u^n+z(\Phi_\tau(u^n)-u^n))-E(u^n)
=\frac12(G(z)-G(0))
+\frac{\varepsilon}{4}(F(z)-F(0)).
\]
Since $P(1)=\delta_E$, the quadratic form of $G$ and integration
by parts give
\begin{align*}
P'(1)
&=\delta_E+\frac12(G'(1)-G(1)+G(0))
+\frac{\varepsilon}{4}(F'(1)-F(1)+F(0))\\
&=\delta_E
+\frac12\|\nabla(\Phi_\tau(u^n)-u^n)\|_2^2
+\frac{\varepsilon}{4}\int_0^1 zF''(z)\,dz.
\end{align*}
Therefore,
\[
\left|
P'(1)-\frac12\|\nabla(\Phi_\tau(u^n)-u^n)\|_2^2
\right|
\le|\delta_E|+\frac{\varepsilon}{4}\int_0^1 z|F''(z)|\,dz
\lesssim\varepsilon\tau^2.
\]

To estimate the effect of normalization, use
$D(1)=\|u^n\|_2^2+\delta_M$ and
$D(1)\ge\|u^n\|_2^2/2$:
\begin{align*}
\left|\frac{\|u^n\|_2^2}{D(1)}-1\right|
&=\frac{|\delta_M|}{D(1)}
\lesssim|\delta_M|,\\
\left|\frac{\|u^n\|_2^4}{D(1)^2}-1\right|
&=\frac{|\delta_M|\bigl(D(1)+\|u^n\|_2^2\bigr)}{D(1)^2}
\lesssim|\delta_M|.
\end{align*}
Moreover, we know that
\begin{align*}
|G(1)-G(0)|
&\le
2\|\nabla u^n\|_2\|\nabla\mathcal N_\tau(u^n)\|_2
+\|\nabla\mathcal N_\tau(u^n)\|_2^2
\lesssim\varepsilon\tau.
\end{align*}
Since $G(0)=\|\nabla u^n\|_2^2$,
\begin{align*}
&\frac{\|u^n\|_2^2G(1)}{2D(1)^2}
+\frac{\varepsilon\|u^n\|_2^4F(1)}{2D(1)^3}
-\frac{\|\nabla u^n\|_2^2}{2\|u^n\|_2^2}\\
&\quad=
\frac{\|u^n\|_2^2}{2D(1)^2}(G(1)-G(0))
+\frac{G(0)}{2\|u^n\|_2^2}
\left(\frac{\|u^n\|_2^4}{D(1)^2}-1\right)
+\frac{\varepsilon\|u^n\|_2^4F(1)}{2D(1)^3}.
\end{align*}
Thus
\begin{align*}
&\left|
\frac{\|u^n\|_2^2G(1)}{2D(1)^2}
+\frac{\varepsilon\|u^n\|_2^4F(1)}{2D(1)^3}
-\frac{\|\nabla u^n\|_2^2}{2\|u^n\|_2^2}
\right|\lesssim
|G(1)-G(0)|+|\delta_M|+\varepsilon|F(1)|
\lesssim\varepsilon.
\end{align*}

Substituting
\[
D'(1)=\|\Phi_\tau(u^n)-u^n\|_2^2+\delta_M,
\qquad
P'(1)=\frac12G'(1)+\frac{\varepsilon}{4}F'(1)
\]
into \eqref{eq:Hderiv1} and using the preceding bounds gives
\begin{align*}
&\left|
H'(1)-P'(1)
+\frac{\|\nabla u^n\|_2^2}{2\|u^n\|_2^2}
\|\Phi_\tau(u^n)-u^n\|_2^2
\right|\\
&\quad\lesssim
|\delta_M|\bigl(|G'(1)|+\varepsilon|F'(1)|+1\bigr)
+\varepsilon\|\Phi_\tau(u^n)-u^n\|_2^2
\lesssim\varepsilon\tau^2.
\end{align*}
Combining this with the estimate for $P'(1)$ yields
\[
\left|
H'(1)-\frac12\left(
\|\nabla(\Phi_\tau(u^n)-u^n)\|_2^2
-\frac{\|\nabla u^n\|_2^2}{\|u^n\|_2^2}
\|\Phi_\tau(u^n)-u^n\|_2^2
\right)
\right|
\lesssim\varepsilon\tau^2.
\]

It remains to compare the centered quadratic increment with
$\Sigma(\tau)$. For the linear increment, Parseval's identity gives
the exact identity
\begin{align*}
&\frac12\left(
\|\nabla[(e^{i\tau\Delta}-I)u^n]\|_2^2
-\frac{\|\nabla u^n\|_2^2}{\|u^n\|_2^2}
\|(e^{i\tau\Delta}-I)u^n\|_2^2
\right)\\
&\quad=
\sum_k(1-\cos(\tau|k|^2))
\left(
|k|^2-\frac{\|\nabla u^n\|_2^2}{\|u^n\|_2^2}
\right)|\widehat u_k^n|^2
=\Sigma(\tau).
\end{align*}
The actual increment differs from the linear one by
$\mathcal N_\tau(u^n)$. Since
\[
\|(e^{i\tau\Delta}-I)u^n\|_{H^s}\lesssim\tau,
\qquad
\|\mathcal N_\tau(u^n)\|_{H^s}\lesssim\varepsilon\tau,
\]
the Cauchy--Schwarz inequality therefore yields
\[
\left|
\frac12\left(
\|\nabla(\Phi_\tau(u^n)-u^n)\|_2^2
-\frac{\|\nabla u^n\|_2^2}{\|u^n\|_2^2}
\|\Phi_\tau(u^n)-u^n\|_2^2
\right)-\Sigma(\tau)
\right|
\lesssim\varepsilon\tau^2.
\]
Combining this estimate with the preceding bound for $H'(1)$ gives
\[
|H'(1)-\Sigma(\tau)|\lesssim\varepsilon\tau^2.
\]
\end{proof}

\subsection{Local existence}
For notational simplicity, we use $\gamma$ and $\lambda$ throughout the remainder of this section to denote the relaxation parameter and the normalization factor associated with the starting value $u^n$, respectively.
We now have enough ingredients to establish the existence of the relaxation parameter $\gamma$ in a neighbourhood of $1$. The main theorem below further estimates the distance between $\gamma$ and $1$. 
\begin{theorem}\label{thm:gamma}
Assume $u^n\in H^{s+\max\{2,\eta\}}$, $s>\max\{1,d/2\}$, $\sigma_N>0$ and $\phi_t(u^n)\in H^{s+\max\{2,\eta\}}$ for $0\leq t\leq\tau$. There are
$\varepsilon_0,\tau_0>0$ and $0<\rho\le1/2$ such that $H$ has exactly
one root in $[1-\rho,1+\rho]$ whenever
$0<\varepsilon\le\varepsilon_0$, $0<\tau\le\tau_0$. This root $\gamma$ satisfies
\begin{equation}
|1-\gamma|\leq C\frac{\varepsilon\tau^{p+1}}{\Sigma(\tau)}\quad\text{and}\quad
|1-\gamma|\leq C\frac{\varepsilon\tau^{p-1}}{\sigma_N}.
\label{eq:gamma2final}
\end{equation}
Moreover, the mass normalization factor $\lambda$ satisfies
\begin{equation}
|1-\lambda|\leq C(1+\sigma_N^{-1})\varepsilon\tau^{p+1}.
\label{eq:lambdaFinal}
\end{equation}
The constant $C>0$ depends on $s,d, \|u^n\|_2$ and $\sup_{0\leq t\leq\tau}\|\phi_t(u^n)\|_{H^{s+\max\{2,\eta\}}}$. The parameters $\varepsilon_0,\tau_0,\rho$ additionally depend on $N$ and $\sigma_N$.
\end{theorem}
\begin{proof}
\emph{Step 1:} We first show the monotonicity of $H$ near $1$.
For the convenience of this proof, let $c\ge1$ dominate the constants in \eqref{eq:H1} and \eqref{eq:Hderiv}, i.e.,
\[
|H(1)|\leq c\varepsilon\tau^{p+1},\quad
|H'(1)-\Sigma(\tau)|\leq c\varepsilon\tau^2,
\qquad
\sup_{|z|\le2}|H''(z)|\leq c\tau^2.
\]
Then for $|z-1|\le\rho$, we have that
\[
|H'(z)-H'(1)|\leq |1-z|\sup_{\xi\in[1,z]}|H''(\xi)|\leq c\rho\tau^2.
\]
Thus, by choosing $\rho=\min\left\{\tfrac12,\frac{\sigma_N}{8c}\right\}$ and $\varepsilon_0 =\min\{1,\frac{\rho\sigma_N}{8c}\}$, we obtain
\begin{align*}
H'(z)&\ge H'(1)-|H'(z)-H'(1)|\\
&\ge\Sigma(\tau)-c(\varepsilon+\rho)\tau^2\\
&\ge\Sigma(\tau)-\tfrac14\sigma_N\tau^2
\ge\tfrac12\Sigma(\tau)\ge\tfrac14\sigma_N\tau^2>0.
\end{align*}

\emph{Step 2:} Next, we show the existence of the root. Since $p\ge1$ and $\tau\le1$,
\[
|H(1)|\le C\varepsilon\tau^{p+1}
\le\tfrac18\rho\sigma_N\tau^2,
\]
Thus,
\begin{align*}
H(1+\rho)&=H(1)+\int_1^{1+\rho}H'(z)\,dz
\ge-\tfrac18\rho\sigma_N\tau^2+\tfrac14\rho\sigma_N\tau^2>0,\\
H(1-\rho)&=H(1)-\int_{1-\rho}^1H'(z)\,dz
\le\tfrac18\rho\sigma_N\tau^2-\tfrac14\rho\sigma_N\tau^2<0.
\end{align*}
Therefore, from the intermediate value theorem, the strict monotonicity gives a unique root $\gamma\in[1-\rho, 1+\rho]$. Moreover,
\[
|H(1)| = |H(\gamma)-H(1)|\geq \tfrac12\Sigma(\tau)|\gamma-1|,
\]
which proves \eqref{eq:gamma2final}.

\emph{Step 3:} Finally, we estimate $|1-\gamma|$. 
\eqref{eq:massidentity} and $\frac12\le\gamma\le\frac32$ imply
\begin{align*}
|D(\gamma)-\|u^n\|_2^2|
&\le |\gamma|\,|\gamma-1|\|\Phi_\tau(u^n)-u^n\|_2^2+|\gamma|\,|\delta_M|\\
&\lesssim |\gamma-1|\tau^2+|\delta_M|.
\end{align*}
The mean value theorem for $x\mapsto\|u^n\|_2x^{-\frac12}$, on
$x\ge\frac12\|u^n\|_2^2$, therefore gives
\[
|\lambda-1|\lesssim |\gamma-1|\tau^2+|\delta_M|
\lesssim(1+\sigma_N^{-1})\varepsilon\tau^{p+1}.
\]
This completes the proof.
\end{proof}

\begin{figure}[htp]
    \centering
    \includegraphics[width=0.5\linewidth]{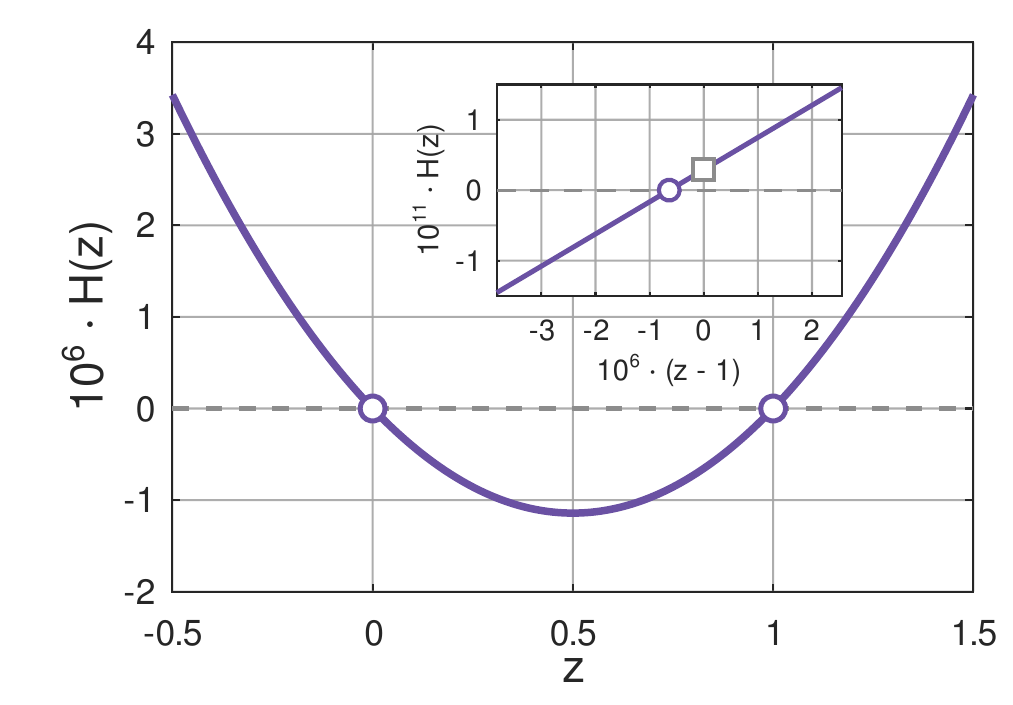}
    \caption{Graph of $H$ in experiments when computing one time step of proposed methods with LRI2 as the baseline integrator.}
    \label{fig:H}
\end{figure}
Figure~\ref{fig:H} illustrates the local uniqueness result of Theorem~\ref{thm:gamma}. On the plotted interval $[-0.5,1.5]$, the graph of $H$ crosses zero at the trivial root $z=0$ and at a single nonzero root $\gamma$. The zoom-in plot further shows that $|1-\gamma|<10^{-6}$ in this test case. Notice that the graph also appears strictly convex throughout the plotted interval. This observation motivates the global uniqueness argument in the next subsection, where we prove that $H$ is strictly convex on an interval containing both roots and is positive outside that interval.

\subsection{Global uniqueness}
%We now prove the stronger statement that $\gamma$ is the only nonzero real relaxation parameter.
Lemma~\ref{lem:incrementupper} provides an upper bound for the increment $\|\Phi_\tau(u^n)-u^n\|_2$ with additional regularity requirements. With the assumption $\sigma_N>0$, we can also establish a uniform lower bound for this increment.
\begin{lemma}\label[lemma]{lem:incrementlower}
Assume $u^n\in H^s$, $s>d/2$, and $\sigma_N>0$. The baseline integrator $\Phi_\tau$ satisfies \eqref{assump:lri}.
For $0<\tau\le\min\{1,\pi/(2N^2)\}$,
\[
\|(e^{i\tau\Delta}-I)u^n\|_2^2
\ge \frac{\sigma_N}{2N^2}\tau^2.
\]
Moreover, for $\varepsilon$ small enough we have
\[
\|\Phi_\tau(u^n)-u^n\|_2^2\ge \frac{\sigma_N}{8N^2}\tau^2.
\]
The smallness threshold for $\varepsilon$ depends on $s,d,N,\sigma_N$ and $\|u^n\|_{H^s}$, and is independent of $\tau$.
\end{lemma}
\begin{proof}
Using the fact that
\[
(a+b)(a-b)^2\le 4N^2(a^2+b^2),\quad0\le a,b\le N^2,
\]
the definition of $\sigma_N$ and Parseval's identity give
\begin{align*}
\sigma_N
&\le \frac{4N^2}{\pi^2\|u^n\|_2^2}
\sum_{|k|,|l|\le N}
(|k|^4+|l|^4)
|\widehat u_k^n|^2|\widehat u_l^n|^2\\
&=\frac{8N^2}{\pi^2\|u^n\|_2^2}
\left(\sum_{|k|\le N}|k|^4|\widehat u_k^n|^2\right)
\left(\sum_{|l|\le N}|\widehat u_l^n|^2\right)\\
&\le \frac{8N^2}{\pi^2}
\sum_{|k|\le N}|k|^4|\widehat u_k^n|^2.
\end{align*}
Since $1-\cos x\ge 2x^2/\pi^2$ for $0\le x\le\pi$,
\begin{align*}
\|(e^{i\tau\Delta}-I)u^n\|_2^2
&=2\sum_k(1-\cos(\tau|k|^2))|\widehat u_k^n|^2\\
&\ge \frac{4\tau^2}{\pi^2}
\sum_{|k|\le N}|k|^4|\widehat u_k^n|^2\\
&\ge \frac{\sigma_N}{2N^2}\tau^2.
\end{align*}
From Assumption~\ref{assump:lri} (i) and (ii), we obtain
\[
\|\mathcal N_\tau(u^n)\|_2\le C\varepsilon\tau.
\]
Finally, by choosing $\varepsilon$ small enough, we obtain
\begin{align*}
\|\Phi_\tau(u^n)-u^n\|_2
&\ge \|(e^{i\tau\Delta}-I)u^n\|_2
-\|\mathcal N_\tau(u^n)\|_2\\
&\ge \left(\sqrt{\frac{\sigma_N}{2N^2}}
-C\varepsilon\right)\tau\\
&\ge \frac12\sqrt{\frac{\sigma_N}{2N^2}}\,\tau,
\end{align*}
which concludes the proof.
\end{proof}
With the aid of the preceding lemma, we derive a strictly positive lower bound for $D$.
\begin{lemma}\label[lemma]{lem:globalmasspositive}
Assume $u^n\in H^{s+\max\{2,\eta\}}$, $s>\max\{1,d/2\}$, $\sigma_N>0$ and $\phi_t(u^n)\in H^{s+\max\{2,\eta\}}$ for $0\leq t\leq\tau$. For sufficiently small $\varepsilon,\tau$,
\begin{equation}
D(z)\ge\tfrac12\left(\|u^n\|_2^2+z^2\|\Phi_\tau(u^n)-u^n\|_2^2\right)
\qquad(z\in\mathbb R).
\label{eq:Mglobalpositive}
\end{equation}
Moreover,
\begin{equation}
\|\Phi_\tau(u^n)-u^n\|_{H^s}\leq C\tau,\qquad
\sup_{z\in\mathbb R}
\left\|\frac{\|u^n\|_2(u^n+z(\Phi_\tau(u^n)-u^n))}{\sqrt{D(z)}}\right\|_{H^s}
\leq C\left(1+\frac{N}{\sqrt{\sigma_N}}\right).
\label{eq:uniformnormalizedHs}
\end{equation}
The constant $C>0$ depends on $s,d$, and $\|u^n\|_{H^{s+\max\{2,\eta\}}}$. The smallness thresholds additionally depend on $\|u^n\|_2,N,\sigma_N$ and $\sup_{0\leq t\leq\tau}\|\phi_t(u^n)\|_{H^{s+\max\{2,\eta\}}}$.
\end{lemma}
\begin{proof}
By Corollary~\ref{lem:defects} and Lemma~\ref{lem:incrementlower},
\[
\frac{|\delta_M|}{\|\Phi_\tau(u^n)-u^n\|_2^2}
\le C\frac{N^2}{\sigma_N}\varepsilon\tau^{p-1}\le\frac12
\]
for sufficiently small $\varepsilon$. The mass identity \eqref{eq:massidentity} then yields
\begin{align*}
D(z)&\ge\|u^n\|_2^2+(z^2-\tfrac32|z|)\|\Phi_\tau(u^n)-u^n\|_2^2\\
&\ge\|u^n\|_2^2+(\tfrac12z^2-\tfrac98)\|\Phi_\tau(u^n)-u^n\|_2^2\\
&\ge\tfrac12\left(\|u^n\|_2^2+z^2\|\Phi_\tau(u^n)-u^n\|_2^2\right),
\end{align*}
after choosing $\tau$ so that $\tfrac98\|\Phi_\tau(u^n)-u^n\|_2^2\le\tfrac12\|u^n\|_2^2$.
Since 
\[
\|\Phi_\tau(u^n)-u^n\|_{H^s}\lesssim\tau,
\]
we obtain that
\begin{align*}
\left\|\frac{\|u^n\|_2(u^n+z(\Phi_\tau(u^n)-u^n))}{\sqrt{D(z)}}\right\|_{H^s}
&\le\sqrt2\|u^n\|_{H^s}
+\sqrt2\|u^n\|_2\frac{\|\Phi_\tau(u^n)-u^n\|_{H^s}}{\|\Phi_\tau(u^n)-u^n\|_2}\\
&\lesssim1+\frac{N}{\sqrt{\sigma_N}}.
\end{align*}
\end{proof}

With all the building blocks established, we show in the following theorem the global uniqueness of $\gamma$.
\begin{theorem}\label{thm:global}
Assume $u^n\in H^{s+\max\{2,\eta\}}$, $s>\max\{1,d/2\}$, $\sigma_N>0$ and $\phi_t(u^n)\in H^{s+\max\{2,\eta\}}$ for $0\leq t\leq\tau$. For sufficiently small $\varepsilon,\tau$, $H$ is defined on $\mathbb R$ and has exactly two real roots: $0$ and the positive root $\gamma$ of Theorem~\ref{thm:gamma}.
The smallness thresholds depend on $s,d,N,\sigma_N,\|u^n\|_2$ and $\sup_{0\leq t\leq\tau}\|\phi_t(u^n)\|_{H^{s+\max\{2,\eta\}}}$.
\end{theorem}
\begin{proof}
Lemma~\ref{lem:globalmasspositive} gives $D(z)>0$ for every $z\in\mathbb R$.
Moreover, Lemmas~\ref{lem:incrementupper} and~\ref{lem:incrementlower}, together with Corollary~\ref{lem:defects}, imply
\[
\|\Phi_\tau(u^n)-u^n\|_{H^s}\lesssim\tau,\qquad
\|\Phi_\tau(u^n)-u^n\|_2\gtrsim\tau,\qquad
|\delta_M|\lesssim\varepsilon\tau^{p+1}.
\]

\emph{Step 1:} We first study the kinetic contribution and define the interval $J_\tau$.
Using \eqref{eq:massidentity} and \eqref{eq:kineticidentity2}, we obtain
\begin{align}
&\frac{\|u^n\|_2^2\|\nabla(u^n+z(\Phi_\tau(u^n)-u^n))\|_2^2}{2D(z)}
-\frac12\|\nabla u^n\|_2^2\nonumber\\
&\quad=\frac{\|u^n\|_2^2z(z-1)}{2D(z)}
\left(\|\nabla(\Phi_\tau(u^n)-u^n)\|_2^2
-\frac{\|\nabla u^n\|_2^2}{\|u^n\|_2^2}\|\Phi_\tau(u^n)-u^n\|_2^2\right)\nonumber\\
&\qquad+\frac{z}{2D(z)}
\left[\|u^n\|_2^2\bigl(\|\nabla\Phi_\tau(u^n)\|_2^2-\|\nabla u^n\|_2^2\bigr)
-\|\nabla u^n\|_2^2\delta_M\right].
\label{eq:kineticFar2}
\end{align}
The comparison of the linear and nonlinear increments in the proof of Lemma~\ref{lem:Hprime} gives
\begin{equation}
\left|\frac12\left(\|\nabla(\Phi_\tau(u^n)-u^n)\|_2^2
-\frac{\|\nabla u^n\|_2^2}{\|u^n\|_2^2}\|\Phi_\tau(u^n)-u^n\|_2^2\right)-\Sigma(\tau)\right|
\lesssim\varepsilon\tau^2.
\label{eq:kineticcoeff2}
\end{equation}
By Proposition~\ref{lem:Sigmalower} and \eqref{eq:kineticcoeff2}, reducing $\varepsilon_0$ gives
\begin{align}
\|\nabla(\Phi_\tau(u^n)-u^n)\|_2^2
-\frac{\|\nabla u^n\|_2^2}{\|u^n\|_2^2}\|\Phi_\tau(u^n)-u^n\|_2^2
\ge2\Sigma(\tau)-C\varepsilon\tau^2
\ge\Sigma(\tau).
\label{eq:positivekinetic2}
\end{align}
For the first quotient in \eqref{eq:kineticFar2}, differentiation of \eqref{eq:massidentity} gives
\begin{align}
\left(\frac{z(z-1)}{D(z)}\right)''
&=\frac{2\|u^n\|_2^2
\left[\|u^n\|_2^2-(3z^2-3z+1)\|\Phi_\tau(u^n)-u^n\|_2^2\right]}{D(z)^3}\nonumber\\
&\quad+\frac{2\delta_M(\|u^n\|_2^2-z^3\|\Phi_\tau(u^n)-u^n\|_2^2)}{D(z)^3}.
\label{eq:quotientsecond2}
\end{align}
We therefore define
\begin{equation}
J_\tau=\left\{z\in\mathbb R:
(3z^2-3z+1)\|\Phi_\tau(u^n)-u^n\|_2^2\le\tfrac23\|u^n\|_2^2\right\}.
\label{eq:Jtau2}
\end{equation}
Since $3z^2-3z+1=3(z-1/2)^2+1/4$, \eqref{eq:Jtau2} and Lemma~\ref{lem:incrementlower} imply
\[
|z|\le\frac12+\sqrt{\frac{2\|u^n\|_2^2}{9\|\Phi_\tau(u^n)-u^n\|_2^2}}
\lesssim\tau^{-1}\qquad(z\in J_\tau).
\]
On the other hand, Lemma~\ref{lem:incrementupper} gives, for $|z|\le2$,
\[
(3z^2-3z+1)\|\Phi_\tau(u^n)-u^n\|_2^2
\le19\|\Phi_\tau(u^n)-u^n\|_2^2
\le\tfrac23\|u^n\|_2^2
\]
for sufficiently small $\tau$. Thus $[-2,2]\subset J_\tau$.

For $z\in J_\tau$, the definition \eqref{eq:Jtau2} also gives
\begin{align*}
z(z-1)\|\Phi_\tau(u^n)-u^n\|_2^2
&=\frac13\left[(3z^2-3z+1)\|\Phi_\tau(u^n)-u^n\|_2^2
-\|\Phi_\tau(u^n)-u^n\|_2^2\right]\\
&\le\frac29\|u^n\|_2^2.
\end{align*}
Combining this bound with \eqref{eq:massidentity}, \eqref{eq:Mglobalpositive} and Corollary~\ref{lem:defects}, we obtain
\begin{align}
\tfrac12\|u^n\|_2^2\le D(z)
&\le\tfrac{11}{9}\|u^n\|_2^2+C\varepsilon\tau^p
\le\tfrac32\|u^n\|_2^2,\nonumber\\
|D'(z)|&\le|\delta_M|+(2|z|+1)\|\Phi_\tau(u^n)-u^n\|_2^2\lesssim\tau,\nonumber\\
D''(z)&=2\|\Phi_\tau(u^n)-u^n\|_2^2\lesssim\tau^2.
\label{eq:DboundsJ2}
\end{align}

\emph{Step 2:} We then show that $H$ is strictly convex on $J_\tau$.
For the second quotient in \eqref{eq:kineticFar2}, differentiation of \eqref{eq:massidentity} gives
\begin{equation}
\left(\frac z{D(z)}\right)''
=\frac{2\left[z^3\|\Phi_\tau(u^n)-u^n\|_2^4
-(3z-1)\|u^n\|_2^2\|\Phi_\tau(u^n)-u^n\|_2^2-\|u^n\|_2^2\delta_M\right]}{D(z)^3}.
\label{eq:quotientlinear2}
\end{equation}
Substituting \eqref{eq:quotientsecond2} and \eqref{eq:quotientlinear2} into the second derivative of \eqref{eq:kineticFar2}, we obtain
\begin{align}
&\frac{\|u^n\|_2^4}{D(z)^3}
\left(\|\nabla(\Phi_\tau(u^n)-u^n)\|_2^2
-\frac{\|\nabla u^n\|_2^2}{\|u^n\|_2^2}\|\Phi_\tau(u^n)-u^n\|_2^2\right)\nonumber\\
&\qquad\cdot\left[\|u^n\|_2^2-(3z^2-3z+1)\|\Phi_\tau(u^n)-u^n\|_2^2\right]
+\mathcal R_1(z).
\label{eq:kineticsecond2}
\end{align}
Here $\mathcal R_1$ contains the contribution from the second term in \eqref{eq:quotientsecond2}
and that from \eqref{eq:quotientlinear2}.

To estimate $\mathcal R_1$, Assumption~\ref{assump:lri} (i), (ii) and $s>1$ first give
$\|\nabla\mathcal N_\tau(u^n)\|_2\lesssim\varepsilon\tau$.
The unitarity of $e^{i\tau\Delta}$ and the Cauchy--Schwarz inequality imply
\begin{align*}
\left|\|\nabla\Phi_\tau(u^n)\|_2^2-\|\nabla u^n\|_2^2\right|
\le2\|\nabla u^n\|_2\|\nabla\mathcal N_\tau(u^n)\|_2
+\|\nabla\mathcal N_\tau(u^n)\|_2^2\lesssim\varepsilon\tau.
\end{align*}
Together with Lemma~\ref{lem:incrementupper} and Corollary~\ref{lem:defects}, this yields
\begin{align}
&\left|\|\nabla(\Phi_\tau(u^n)-u^n)\|_2^2
-\frac{\|\nabla u^n\|_2^2}{\|u^n\|_2^2}\|\Phi_\tau(u^n)-u^n\|_2^2\right|
\lesssim\tau^2,\nonumber\\
&\left|\|u^n\|_2^2\bigl(\|\nabla\Phi_\tau(u^n)\|_2^2-\|\nabla u^n\|_2^2\bigr)
-\|\nabla u^n\|_2^2\delta_M\right|
\lesssim\varepsilon\tau.
\label{eq:kineticbounds2}
\end{align}
Using \eqref{eq:quotientsecond2}, \eqref{eq:quotientlinear2}, \eqref{eq:kineticbounds2} and \eqref{eq:DboundsJ2}, we obtain
\begin{align*}
|\mathcal R_1(z)|
&\lesssim\varepsilon\tau^4(1+\tau^2|z|^3)
+\varepsilon\tau\bigl(\tau^4|z|^3+\tau^2(1+|z|)+\varepsilon\tau^2\bigr)\\
&\lesssim\varepsilon\tau^4(1+\tau^{-1})
+\varepsilon\tau\bigl(\tau+\tau^2+\tau+\varepsilon\tau^2\bigr)
\lesssim\varepsilon\tau^2,
\end{align*}
where we used $|z|\lesssim\tau^{-1}$ from \eqref{eq:Jtau2}.
Moreover, \eqref{eq:Jtau2} and \eqref{eq:DboundsJ2} imply
\[
\|u^n\|_2^2-(3z^2-3z+1)\|\Phi_\tau(u^n)-u^n\|_2^2\ge\tfrac13\|u^n\|_2^2,
\qquad D(z)^3\le\tfrac{27}{8}\|u^n\|_2^6.
\]
Therefore, \eqref{eq:kineticsecond2} and \eqref{eq:positivekinetic2} show that the second derivative of the kinetic contribution is bounded below by
\[
\frac8{81}\Sigma(\tau)-C\varepsilon\tau^2.
\]

For the quartic contribution, define
\[
F(z)=\|u^n+z(\Phi_\tau(u^n)-u^n)\|_4^4.
\]
By Lemma~\ref{lem:incrementupper} and \eqref{eq:Jtau2}, for $z\in J_\tau$ we have
\[
\|u^n+z(\Phi_\tau(u^n)-u^n)\|_{H^s}
\le\|u^n\|_{H^s}+|z|\|\Phi_\tau(u^n)-u^n\|_{H^s}
\lesssim1.
\]
The Sobolev embedding $H^s\hookrightarrow L^\infty$ and H\"older's inequality give
\begin{align*}
|F(z)|&\lesssim1,\\
|F'(z)|&\le4\|u^n+z(\Phi_\tau(u^n)-u^n)\|_\infty^2
\|u^n+z(\Phi_\tau(u^n)-u^n)\|_2\|\Phi_\tau(u^n)-u^n\|_2
\lesssim\tau,\\
|F''(z)|&\le12\|u^n+z(\Phi_\tau(u^n)-u^n)\|_\infty^2
\|\Phi_\tau(u^n)-u^n\|_2^2
\lesssim\tau^2.
\end{align*}
Using \eqref{eq:DboundsJ2}, we then obtain
\begin{align*}
\left(\frac F{D^2}\right)''
&=\frac{F''}{D^2}-\frac{4D'F'}{D^3}
-\frac{2D''F}{D^3}+\frac{6(D')^2F}{D^4},\\
\left|\frac{\varepsilon\|u^n\|_2^4}{4}
\left(\frac F{D^2}\right)''\right|
&\lesssim\varepsilon(\tau^2+\tau\cdot\tau+\tau^2+\tau^2)
\lesssim\varepsilon\tau^2.
\end{align*}
Combining this estimate with the kinetic lower bound from \eqref{eq:kineticsecond2} in \eqref{eq:Hexpanded1} gives
\begin{equation}
H''(z)\ge\frac8{81}\Sigma(\tau)-C\varepsilon\tau^2
\ge\frac1{12}\Sigma(\tau)>0\qquad(z\in J_\tau),
\label{eq:HconvexJ2}
\end{equation}
after choosing $C\varepsilon_0\le5\sigma_N/648$, where Proposition~\ref{lem:Sigmalower} is used in the last inequality.

\emph{Step 3:} Next, we show that $H$ is positive outside $J_\tau$.
If $z\notin J_\tau$, then $|z|>2$ since $[-2,2]\subset J_\tau$.
By \eqref{eq:Jtau2} and Lemma~\ref{lem:incrementupper},
\[
z(z-1)\|\Phi_\tau(u^n)-u^n\|_2^2
>\frac29\|u^n\|_2^2-\frac13\|\Phi_\tau(u^n)-u^n\|_2^2
\ge\frac19\|u^n\|_2^2
\]
for sufficiently small $\tau$.
Since $z(z-1)\le2z^2$ when $|z|>2$, Lemma~\ref{lem:incrementupper} further gives
\[
\frac19\|u^n\|_2^2
<z(z-1)\|\Phi_\tau(u^n)-u^n\|_2^2
\lesssim z^2\tau^2.
\]
Thus $|z|\gtrsim\tau^{-1}$ and $|z-1|\ge|z|/2$.
Using Corollary~\ref{lem:defects} and Lemma~\ref{lem:incrementlower}, we obtain
\begin{align*}
\frac{|z\delta_M|}{z(z-1)\|\Phi_\tau(u^n)-u^n\|_2^2}
=\frac{|\delta_M|}{|z-1|\|\Phi_\tau(u^n)-u^n\|_2^2}\le C\frac{\varepsilon\tau^{p+1}}{\tau^{-1}\tau^2}
=C\varepsilon\tau^p\le\frac12.
\end{align*}
The mass identity \eqref{eq:massidentity} therefore implies
\begin{align}
D(z)
&\le\|u^n\|_2^2+\tfrac32z(z-1)\|\Phi_\tau(u^n)-u^n\|_2^2\nonumber\\
&\le\tfrac{21}{2}z(z-1)\|\Phi_\tau(u^n)-u^n\|_2^2
\le15z(z-1)\|\Phi_\tau(u^n)-u^n\|_2^2.
\label{eq:DoutsideJ2}
\end{align}
By \eqref{eq:DoutsideJ2} and \eqref{eq:positivekinetic2}, the first term on the right hand side of \eqref{eq:kineticFar2} is at least
\begin{align*}
&\frac{\|u^n\|_2^2}{30\|\Phi_\tau(u^n)-u^n\|_2^2}
\left(\|\nabla(\Phi_\tau(u^n)-u^n)\|_2^2
-\frac{\|\nabla u^n\|_2^2}{\|u^n\|_2^2}\|\Phi_\tau(u^n)-u^n\|_2^2\right)\\
&\quad\ge\frac{\|u^n\|_2^2\Sigma(\tau)}
{30\|\Phi_\tau(u^n)-u^n\|_2^2}
\ge2c_0>0.
\end{align*}
Here $c_0>0$ is independent of $\varepsilon,\tau$, since Lemma~\ref{lem:incrementupper}
gives $\|\Phi_\tau(u^n)-u^n\|_2^2\le C\tau^2$ and Proposition~\ref{lem:Sigmalower}
gives $\Sigma(\tau)\ge\sigma_N\tau^2/2$.

For the second term in \eqref{eq:kineticFar2}, \eqref{eq:Mglobalpositive} and Lemma~\ref{lem:incrementlower} give
\[
\frac{|z|}{D(z)}
\le\frac{2|z|}{\|u^n\|_2^2+z^2\|\Phi_\tau(u^n)-u^n\|_2^2}
\le\frac1{\|u^n\|_2\|\Phi_\tau(u^n)-u^n\|_2}
\lesssim\tau^{-1}.
\]
Together with \eqref{eq:kineticbounds2}, this bounds its absolute value by
$C\varepsilon\tau\cdot\tau^{-1}\lesssim\varepsilon$.
Finally, \eqref{eq:uniformnormalizedHs} and the Sobolev embedding $H^s\hookrightarrow L^\infty$ imply
\begin{align*}
\left|\frac{\varepsilon\|u^n\|_2^4F(z)}{4D(z)^2}
-\frac\varepsilon4\|u^n\|_4^4\right|\le\frac\varepsilon4\left(
\left\|\frac{\|u^n\|_2[u^n+z(\Phi_\tau(u^n)-u^n)]}{\sqrt{D(z)}}\right\|_4^4
+\|u^n\|_4^4\right)
\lesssim\varepsilon.
\end{align*}
Combining these estimates in \eqref{eq:kineticFar2} and \eqref{eq:Hexpanded1}, we obtain
\begin{equation}
H(z)\ge2c_0-C\varepsilon\ge c_0>0\qquad(z\notin J_\tau),
\label{eq:Hpositiveoutside2}
\end{equation}
after reducing $\varepsilon_0$.

Since the roots $0$ and $\gamma$ from Corollary~\ref{lem:H1} and Theorem~\ref{thm:gamma} belong to $J_\tau$,
\eqref{eq:HconvexJ2} and \eqref{eq:Hpositiveoutside2} exclude every other real root of $H$.
\end{proof}

\section{Convergence analysis}\label{sec:convergence_analysis}
We are now in the position to study the convergence of the proposed scheme. We begin with a local error estimate for a general input $v$ and prescribed parameters $\gamma$ and $\lambda$. To this end, we define for any $\lambda,\gamma>0$, and $s>d/2$,
\begin{align*}
    \Psi_{\lambda, \gamma} : H^s&\rightarrow H^s\\
    v&\mapsto \lambda(v+\gamma(\Phi_\tau(v) - v)).
\end{align*}
By Assumption~\ref{assump:lri}, $\Phi_\tau$ maps $H^s$ into itself, so $\Psi_{\lambda,\gamma}$ is well defined.
\begin{proposition}\label[proposition]{prop:local_error}
Let $s>d/2$, $0<\varepsilon,\tau\leq1$ and $v\in H^{s+\max\{2,\eta\}}$. Suppose that $\Phi_\tau$ satisfies Assumption~\ref{assump:lri} and that
\[
\phi_t(v)\in H^{s+\max\{2,\eta\}},\qquad 0\leq t\leq3\tau/2.
\]
Then, for any $\gamma\in[\frac12,\frac32]$ and $\lambda>0$,
\begin{equation}
\left\|\Psi_{\lambda,\gamma}(v)-\phi_{\gamma\tau}(v)\right\|_{H^s}
\leq C\left(\varepsilon\tau^{p+1}+|1-\gamma|\tau+|1-\lambda|\right),
\label{eq:localerror2}
\end{equation}
where $C>0$ depends on $s,d$, and $\sup_{0\leq t\leq3\tau/2}\|\phi_t(v)\|_{H^{s+\max\{2,\eta\}}}$.
Moreover, if
\[
\phi_t(v)\in H^{s+\max\{4,\eta\}},\qquad 0\leq t\leq3\tau/2,
\]
then
\begin{equation}
\left\|\Psi_{\lambda,\gamma}(v)-\phi_{\gamma\tau}(v)\right\|_{H^s}
\leq C\left(\varepsilon\tau^{p+1}+|1-\gamma|\tau^2+|1-\lambda|\right),
\label{eq:localerror2smooth}
\end{equation}
where $C>0$ depends on $s,d$, and $\sup_{0\leq t\leq3\tau/2}\|\phi_t(v)\|_{H^{s+\max\{4,\eta\}}}$. Both constants are independent of $\gamma$ and $\lambda$.
\end{proposition}
\begin{proof}
    We first split the local error into several terms.
    \[
    \Psi_{\lambda, \gamma}(v)-\phi_{\gamma\tau}(v) = 
    (\lambda-1)\bigl[v+\gamma(\Phi_\tau(v)-v)\bigr]
    +\gamma\bigl[\Phi_\tau(v)-\phi_\tau(v)\bigr]
    +(\gamma-1)\bigl[\phi_\tau(v)-v\bigr]
    -\bigl[\phi_{\gamma\tau}(v)-\phi_\tau(v)\bigr].
    \]
    The first two terms can be estimated directly. Indeed,
    \begin{equation}\label{eq:local_error12}
    \|(\lambda-1)\bigl[v+\gamma(\Phi_\tau(v)-v)\bigr]\|_{H^s}\lesssim |1-\lambda|, \quad \|\gamma\bigl[\Phi_\tau(v)-\phi_\tau(v)\bigr]\|_{H^s}\lesssim\varepsilon\tau^{p+1}.
    \end{equation}
    We know that
    \begin{align*}
    \partial_t\phi_t(v)&=i\Delta\phi_t(v) - i\varepsilon|\phi_t(v)|^2\phi_t(v),\\
        \partial_t^2\phi_t(v)&=i\Delta\partial_t\phi_t(v)
    -i\varepsilon\left(2|\phi_t(v)|^2\partial_t\phi_t(v)
    +\phi_t(v)^2\overline{\partial_t\phi_t(v)}\right).
    \end{align*}
    Thus, 
    \begin{align*}
    \|\partial_t\phi_t(v)\|_{H^s}&\lesssim\|\phi_t(v)\|_{H^{s+2}}+\varepsilon\|\phi_t(v)\|^3_{H^s},\\
    \|\partial_t^2\phi_t(v)\|_{H^s}&\lesssim \|\phi_t(v)\|_{H^{s+4}} + \varepsilon\|\phi_t(v)\|^3_{H^{s+2}}.
    \end{align*}
    If $v$ only satisfies $\phi_t(v)\in H^{s+\max\{2,\eta\}}$ for $0\leq t\leq3\tau/2$, the last two terms can be estimated as
    \[
    \|(\gamma-1)\bigl[\phi_\tau(v)-v\bigr]\|_{H^s}\leq|1-\gamma|\|\int_0^\tau\partial_t\phi_t(v)\dd t\|_{H^s}\leq |1-\gamma|\int_0^\tau\|\partial_t\phi_t(v)\|_{H^s}\dd t\lesssim|1-\gamma|\tau
    \]
    and
    \[
    \|\phi_{\gamma\tau}(v)-\phi_\tau(v)\|_{H^s}=\|\int_\tau^{\gamma\tau}\partial_t\phi_t(v)\dd t\|_{H^s}\leq \left|\int_\tau^{\gamma\tau}\|\partial_t\phi_t(v)\|_{H^s}\dd t\right|\lesssim|1-\gamma|\tau.
    \]
    Combining these estimates with \eqref{eq:local_error12}, we arrive at the local error estimate \eqref{eq:localerror2}.

    If $v$ further satisfies $\phi_t(v)\in H^{s+\max\{4,\eta\}}$ for $0\leq t\leq3\tau/2$, the last two terms can be estimated as
    \begin{align*}
        &\|(\gamma-1)\bigl[\phi_\tau(v)-v\bigr] - \phi_{\gamma\tau}(v)+\phi_\tau(v)\|_{H^s} \\
        &= \|(\gamma-1)\int_0^\tau\left(\partial_t\phi_t(v) - \partial_t\phi_\tau(v)\right)\dd t - \int_\tau^{\gamma\tau}\left(\partial_t\phi_t(v) - \partial_t\phi_\tau(v)\right)\dd t\|_{H^s}\\
        &\leq |1-\gamma|\int_0^\tau\|\partial_t\phi_t(v) - \partial_t\phi_\tau(v)\|_{H^s}\dd t + \left|\int_\tau^{\gamma\tau}\|\partial_t\phi_t(v) - \partial_t\phi_\tau(v)\|_{H^s}\dd t\right|.
    \end{align*}
    It remains to estimate $\|\partial_t\phi_t(v) - \partial_t\phi_\tau(v)\|_{H^s}$. In fact,
    \[
    \|\partial_t\phi_t(v) - \partial_t\phi_\tau(v)\|_{H^s} = \|\int_\tau^t\partial_t^2\phi_\sigma(v)\dd \sigma\|_{H^s}\leq \sup_{0\leq\sigma\leq3\tau/2}\|\partial_t^2\phi_\sigma(v)\|_{H^s}|t-\tau|\lesssim|t-\tau|.
    \]
    Hence, we know that 
    \[
    \|(\gamma-1)\bigl[\phi_\tau(v)-v\bigr] - \phi_{\gamma\tau}(v)+\phi_\tau(v)\|_{H^s}\lesssim|1-\gamma|\tau^2.
    \]
    Combining these estimates with \eqref{eq:local_error12}, we arrive at the local error estimate \eqref{eq:localerror2smooth}.
\end{proof}
The following corollary is the direct consequence of Proposition~\ref{prop:local_error} and gives a local error estimate for one step of the proposed method starting from $u^n$.
\begin{corollary}[Local error]\label{cor:local_error}
Suppose that $u^n\in H^{s+\max\{2,\eta\}}$, $s>\max\{1,d/2\}$, and that $\Phi_\tau$ satisfies
Assumption~\ref{assump:lri}. Assume that $\sigma_N(u^n)>0$ for some
$N\in\mathbb N$, and that $\phi_t(u^n)\in H^{s+\max\{2,\eta\}}$ for $0\le t\le 3\tau/2$.
For sufficiently small $\varepsilon,\tau>0$, let $\gamma_n$ be the
nonzero root of $H$ near $1$ and let $\lambda_n$ be the corresponding
mass normalization factor. Then
\[
\|u^{n+1}-\phi_{\gamma_n\tau}(u^n)\|_{H^s}
\le C\varepsilon\tau^p,
\]
where  $C$ depends on $s,d,\|u^n\|_2$,$\sigma_N$, and $\sup_{0\leq t\leq3\tau/2}\|\phi_t(u^n)\|_{H^{s+\max\{2,\eta\}}}$.
If, in addition, $\phi_t(u^n)\in H^{s+\max\{4,\eta\}}$ for
$0\le t\le 3\tau/2$, then
\[
\|u^{n+1}-\phi_{\gamma_n\tau}(u^n)\|_{H^s}
\le C\varepsilon\tau^{p+1}.
\]
where $C$ depends on $s,d,\|u^n\|_2$,
$\sigma_N$, and $\sup_{0\leq t\leq3\tau/2}\|\phi_t(u^n)\|_{H^{s+\max\{4,\eta\}}}$.
\end{corollary}

\begin{proof}
Theorem~\ref{thm:gamma} gives $\gamma_n\in[\frac12,\frac32]$ and
\[
|1-\gamma_n|\lesssim\varepsilon\tau^{p-1},
\qquad
|1-\lambda_n|\lesssim\varepsilon\tau^{p+1}.
\]
Applying Proposition~\ref{prop:local_error} with $v=u^n$ and the
computed parameters $\gamma_n$ and $\lambda_n$ gives the two bounds.
\end{proof}

Using the local error estimates in Corollary~\ref{cor:local_error}, we now derive a global error bound for scheme~\eqref{eq:scheme2} with baseline integrators of order at least two, comparing the numerical iterates with the exact solution $u(t)=\phi_t(u^0)$ at the corresponding time nodes.
\begin{theorem}[Global error]\label{thm:globalerror2}
Let $s>\max\{1,d/2\}$ and suppose that the baseline scheme $\Phi_\tau$ satisfies Assumption~\ref{assump:lri} and $p\ge2$. Assume that, for some fixed $N\in\mathbb N$ and $\sigma_*>0$ independent of $\varepsilon$ and $\tau$,
\[
\inf_{t_n\le T}\sigma_N(u^n)\ge\sigma_*.
\]
Suppose that $u^n\in H^{s+\max\{2,\eta\}}$ for every $n$ with $t_n\le T$ and that the exact solution of \eqref{eq:nlse} satisfies
\[
u(t)\in H^{s+\max\{2,\eta\}},\qquad 0\leq t\leq T+2\tau.
\]
Then, for sufficiently small $\varepsilon,\tau>0$, we have that
\begin{equation}
\begin{aligned}
\|u^n-u(t_n)\|_{H^s}\leq C T e^{2LT}\left(\varepsilon\tau^p
+\max_{0\le k<n}|1-\gamma_k|
+\tau^{-1}\max_{0\le k<n}|1-\lambda_k|\right),
\end{aligned}
\label{eq:globalerror2low}
\end{equation}
where $C,L>0$ depend on $s,d,N,\sigma_*,\|u^0\|_2$, and $\sup_{0\leq t\leq T+2\tau}\|u(t)\|_{H^{s+\max\{2,\eta\}}}$.
Moreover, if the exact solution of \eqref{eq:nlse} satisfies
\[
u(t)\in H^{s+\max\{4,\eta\}},\qquad 0\leq t\leq T+2\tau,
\]
then
\begin{equation}
\begin{aligned}
\|u^n-u(t_n)\|_{H^s}\leq C T e^{2LT}\left(\varepsilon\tau^p
+\max_{0\le k<n}|1-\gamma_k|\tau
+\tau^{-1}\max_{0\le k<n}|1-\lambda_k|\right),
\end{aligned}
\label{eq:globalerror2smooth}
\end{equation}
where $C,L>0$ depend on $s,d,N,\sigma_*,\|u^0\|_2$ and $\sup_{0\leq t\leq T+2\tau}\|u(t)\|_{H^{s+\max\{4,\eta\}}}$. The smallness thresholds additionally depend on $T$ and the initial datum $u^0$.
\end{theorem}
\begin{proof}
    Set $R=\sup_{0\leq t\leq T+2\tau}\|u(t)\|_{H^{s+\max\{2,\eta\}}}$.
    We first work under the bound
    $\|u^n\|_{H^{s+\max\{2,\eta\}}}\le R+1$,
    which will be verified below. For sufficiently small $\tau$, local well-posedness gives
    \[
    \sup_{0\leq t\leq\tau}\|\phi_t(u^n)\|_{H^{s+\max\{2,\eta\}}}\leq2(R+1),
    \]
    uniformly in $n$ and $0<\varepsilon\leq1$.
    Thus, applying Theorem~\ref{thm:gamma} and the lower bound $\sigma_N(u^n)\ge\sigma_*$, we obtain the estimates
    \[
    \tfrac12\le\gamma_n\le\tfrac32,\qquad
    |1-\gamma_n|\lesssim\varepsilon\tau^{p-1},\qquad
    |1-\lambda_n|\lesssim\varepsilon\tau^{p+1},
    \]
    uniformly for $t_n\leq T$.

    Assumption~\ref{assump:lri} (ii) implies that for any $v,v'\in H^s$ with $\|v\|_{H^s},\|v'\|_{H^s}\le R+1$,
    \[
    \|\Phi_\tau(v)-\Phi_\tau(v')\|_{H^s}
    \le(1+C\varepsilon\tau)\|v-v'\|_{H^s}.
    \]
    To estimate the difference between the two outputs, we apply the relaxed scheme to both inputs with the same relaxation parameter $\gamma_n$ and normalization factor $\lambda_n$. The stability estimate then gives
    \begin{align*}
    \|\Psi_{\lambda_n, \gamma_n}(v)- \Psi_{\lambda_n, \gamma_n}(v')\|_{H^s}
    &=\|\lambda_n[v+\gamma_n(\Phi_\tau(v)-v)]
    -\lambda_n[v'+\gamma_n(\Phi_\tau(v')-v')]\|_{H^s}\\
    &=\|\lambda_n(1-\gamma_n)(v-v')
    +\gamma_n\lambda_n(\Phi_\tau(v)-\Phi_\tau(v'))\|_{H^s}\\
    &\le\lambda_n\bigl(|1-\gamma_n|
    +\gamma_n(1+C\varepsilon\tau)\bigr)\|v-v'\|_{H^s}.
    \end{align*}
    The multiplier satisfies
    \begin{align*}
    \lambda_n\bigl(|1-\gamma_n|
    +\gamma_n(1+C\varepsilon\tau)\bigr)
    &\le(1+|1-\lambda_n|)
    \bigl(1+2|1-\gamma_n|+\tfrac32C\varepsilon\tau\bigr)\\
    &\le\exp\left(|1-\lambda_n|+2|1-\gamma_n|
    +\tfrac32C\varepsilon\tau\right)\\
    &\le e^{L\tau}.
    \end{align*}
    Here, $L>0$ is independent of $\varepsilon,\tau$. By adapting $L$, the stability estimate
    \[
    \|\Psi_{\lambda_n, \gamma_n}(v)- \Psi_{\lambda_n, \gamma_n}(v')\|_{H^{s+\max\{2,\eta\}}}\leq e^{L\tau}\|v-v'\|_{H^{s+\max\{2,\eta\}}}
    \]
    also holds for $\|v\|_{H^{s+\max\{2,\eta\}}},\|v'\|_{H^{s+\max\{2,\eta\}}}\le R+1$.

    We now verify the bound $\|u^n\|_{H^{s+\max\{2,\eta\}}}\le R+1$. By continuous dependence of the exact flow in $H^{s+\max\{2,\eta\}}$ and persistence of regularity, for any sufficiently small $\delta>0$ we can choose smooth initial data $w_0$ such that $w(t)=\phi_t(w_0)$ satisfies
    \[
    \sup_{0\le t\le T+2\tau}
    \|w(t)-u(t)\|_{H^{s+\max\{2,\eta\}}}\le\delta.
    \]
    In particular, $\|w(t)\|_{H^{s+\max\{2,\eta\}}}\le R+\delta$.
    Applying Proposition~\ref{prop:local_error} with Sobolev regularity $s+\max\{2,\eta\}$ to $w$ gives
    \begin{align*}
    \|\Psi_{\lambda_n,\gamma_n}(w(t_n))
    -w(t_{n+1})\|_{H^{s+\max\{2,\eta\}}}\le C_w\left(\varepsilon\tau^{p+1}
    +|1-\gamma_n|\tau+|1-\lambda_n|\right)
    \le C_w\tau^2.
    \end{align*}
    Here, $C_w$ depends on $w$, but is independent of $n,\varepsilon,\tau$. Using the same parameters on $u^n$ and $w(t_n)$, the stability estimate in $H^{s+\max\{2,\eta\}}$ yields
    \[
    \|u^{n+1}-w(t_{n+1})\|_{H^{s+\max\{2,\eta\}}}
    \le e^{L\tau}\|u^n-w(t_n)\|_{H^{s+\max\{2,\eta\}}}
    +C_w\tau^2.
    \]
    Since $n\tau\le2t_n\le2T$, iteration gives
    \[
    \|u^n-w(t_n)\|_{H^{s+\max\{2,\eta\}}}
    \le e^{2LT}\left(\delta+2TC_w\tau\right).
    \]
    For sufficiently small $\delta$ and $\tau$, we obtain
    \[
    \|u^n-u(t_n)\|_{H^{s+\max\{2,\eta\}}}
    \le e^{2LT}\left(\delta+2TC_w\tau\right)+\delta
    \le\tfrac12.
    \]
    Hence $\|u^n\|_{H^{s+\max\{2,\eta\}}}\le R+\tfrac12$. This closes the bootstrap by induction, and the above parameter and stability estimates hold uniformly for all $t_n\le T$.

    Therefore, since
    \begin{align*}
    u^{n+1}-u(t_{n+1})
    &=\lambda_n\bigl[u^n+\gamma_n(\Phi_\tau(u^n)-u^n)\bigr]\\
    &\quad-\lambda_n\bigl[u(t_n)
    +\gamma_n(\Phi_\tau(u(t_n))-u(t_n))\bigr]\\
    &\quad+\bigl(\Psi_{\lambda_n,\gamma_n}(u(t_n))
    -\phi_{\gamma_n\tau}(u(t_n))\bigr),
    \end{align*}
    Proposition~\ref{prop:local_error} gives
    \begin{equation}
    \begin{aligned}
    \|u^{n+1}-u(t_{n+1})\|_{H^s}
    \le e^{L\tau}\|u^n-u(t_n)\|_{H^s}
    +C\left(\varepsilon\tau^{p+1}
    +|1-\gamma_n|\tau+|1-\lambda_n|\right).
    \end{aligned}
    \label{eq:errorrecurrence2}
    \end{equation}
    Under the stronger assumption
    $u(t)\in H^{s+\max\{4,\eta\}}$ for $0\leq t\leq T+2\tau$,
    the sharper recurrence is
    \begin{equation}
    \begin{aligned}
    \|u^{n+1}-u(t_{n+1})\|_{H^s}
    \le e^{L\tau}\|u^n-u(t_n)\|_{H^s}
    +C\left(\varepsilon\tau^{p+1}
    +|1-\gamma_n|\tau^2+|1-\lambda_n|\right).
    \end{aligned}
    \label{eq:errorrecurrence2smooth}
    \end{equation}
    Therefore, iterating \eqref{eq:errorrecurrence2} gives
    \begin{align*}
    \|u^n-u(t_n)\|_{H^s}
    &\le C\sum_{k=0}^{n-1}e^{L(n-1-k)\tau}
    \left(\varepsilon\tau^{p+1}
    +|1-\gamma_k|\tau+|1-\lambda_k|\right).
    \end{align*}
    The bound $\gamma_k\ge1/2$ (from Theorem~\ref{thm:gamma}) implies that
    \[
    n\tau\le2t_n\le2T,\qquad
    e^{L(n-1-k)\tau}\le e^{2LT}.
    \]
    Thus, we arrive at
    \begin{align*}
    \|u^n-u(t_n)\|_{H^s}
    \lesssim T e^{2LT}\left(\varepsilon\tau^p
    +\max_{k<n}|1-\gamma_k|
    +\tau^{-1}\max_{k<n}|1-\lambda_k|\right).
    \end{align*}
    Iterating \eqref{eq:errorrecurrence2smooth} in the same
    way gives
    \begin{align*}
    \|u^n-u(t_n)\|_{H^s}
    \lesssim T e^{2LT}\left(\varepsilon\tau^p
    +\tau\max_{k<n}|1-\gamma_k|
    +\tau^{-1}\max_{k<n}|1-\lambda_k|\right).
    \end{align*}
    All implicit constants are uniform in $n,\varepsilon,\tau$.
\end{proof}

\section{Numerical examples}\label{sec:numerics}
In this section, we assess the $H^1$ accuracy, conservation properties, and computational cost of the mass-normalized time-relaxation scheme \eqref{eq:scheme2} for the one-dimensional cubic NLSE \eqref{eq:nlse}. Our main method, TRLRI2, uses the second-order LRI2 \eqref{eq:nlse_lri2} as its baseline. We compare TRLRI2 with the uncorrected LRI2 and Strang splitting methods, the time-relaxed Strang splitting method (TRStrang), and Besse's relaxation Crank--Nicolson method (relaxCN)~\cite{CN_relax}.

\textbf{Comparison methods.}
The exact flows of the linear and nonlinear subproblems are
\[
\varphi_\tau^{\mathcal L}(u)=\e^{i\tau\Delta}u,
\qquad
\varphi_\tau^{\mathcal N}(u)=\e^{-i\varepsilon\tau|u|^2}u.
\]
For TRStrang, we apply the mass normalization and energy correction of \eqref{eq:scheme2} to a Strang splitting step,
\begin{align*}
\tilde u^{n+1}
&=\varphi_{\tau/2}^{\mathcal L}\circ\varphi_\tau^{\mathcal N}\circ\varphi_{\tau/2}^{\mathcal L}(u^n),\\
u^{n+1}
&=\frac{\|u^n\|_2\bigl[u^n+\gamma_n(\tilde u^{n+1}-u^n)\bigr]}
{\|u^n+\gamma_n(\tilde u^{n+1}-u^n)\|_2}.
\end{align*}
We choose $\gamma_n$ as the nonzero real root near $1$ of the equation $\E(u^{n+1})=\E(u^n)$ and advance time by $\gamma_n\tau$, as in \eqref{eq:scheme2}.

The relaxCN scheme introduces a real-valued auxiliary variable $\phi^{n+1/2}$ that approximates the density $|u|^2$ at time $(n+1/2)\tau$. The scheme is given by
\begin{equation}\label{eq:relaxCN}
\begin{cases}
\displaystyle
\frac{\phi^{n+1/2}+\phi^{n-1/2}}{2}=|u^n|^2,\\[0.6em]
\displaystyle
i\frac{u^{n+1}-u^n}{\tau}
+\Delta\left(\frac{u^{n+1}+u^n}{2}\right)
=\varepsilon\phi^{n+1/2}\frac{u^{n+1}+u^n}{2}.
\end{cases}
\end{equation}
The first equation determines $\phi^{n+1/2}$ explicitly, and the second is then linear in $u^{n+1}$. In our experiments, we set $u^0=u_0$ and $\phi^{-1/2}=|u^0|^2$. For smooth solutions, relaxCN is expected to be second order in time. The second-order $H^s$ error estimate in~\cite{Besse2021} assumes initial data in $H^{s+4}$ with $s>1/2$, temporal smoothness, uniform bounds on the numerical solution, and a second-order approximation of $|u(-\tau/2)|^2$ in the initialization of $\phi^{-1/2}$.

In addition to mass, relaxCN preserves the discrete energy
\[
\frac12\|\nabla u^n\|_2^2
+\frac{\varepsilon}{4}\int_{\mathbb T}
\left(2\phi^{n-1/2}|u^n|^2-(\phi^{n-1/2})^2\right)\,dx
=\E(u^n)-\frac{\varepsilon}{4}\bigl\||u^n|^2-\phi^{n-1/2}\bigr\|_2^2.
\]
With the initialization used here, this invariant equals $\E(u^0)$ and the physical energy $\E(u^n)$ can vary with $n$.

\textbf{Computation of the relaxation parameter.}
For both TRLRI2 and TRStrang, the energy equation \eqref{eq:H} can be reduced to a polynomial equation of degree at most three.
\begin{proposition}\label[proposition]{prop:polynomial}
Suppose that $D(z)>0$ for every $z\in\mathbb R$, and define
\begin{align}\label{eq:P}
P(z)
&=2\|u^n\|_2^2D(z)
\bigl\|\nabla[u^n+z(\Phi_\tau(u^n)-u^n)]\bigr\|_2^2\nonumber\\
&\quad+\varepsilon\|u^n\|_2^4
\|u^n+z(\Phi_\tau(u^n)-u^n)\|_4^4
-4\E(u^n)D(z)^2.
\end{align}
Then $P(z)=zQ(z)$ for a real polynomial $Q$ of degree at most three. For every real $z\ne0$,
\[
H(z)=0\quad\Longleftrightarrow\quad Q(z)=0.
\]
\end{proposition}
\begin{proof}
Multiplying \eqref{eq:Hexpanded1} by $4D(z)^2$ gives $P(z)=4D(z)^2H(z)$. As functions of the real variable $z$, both $D(z)$ and
$\|\nabla[u^n+z(\Phi_\tau(u^n)-u^n)]\|_2^2$
are real polynomials of degree at most two, and
$\|u^n+z(\Phi_\tau(u^n)-u^n)\|_4^4$
is a real polynomial of degree at most four. Thus $P$ is a real polynomial of degree at most four. Since $H(0)=0$, we have $P(0)=0$, so $P(z)=zQ(z)$ with $\deg Q\le3$. The equivalence follows from $D(z)>0$ and $z\ne0$.
\end{proof}

Lemma~\ref{lem:globalmasspositive} guarantees the positivity of $D$ under its stated hypotheses. In the experiments, we compute the roots of $Q$ with MATLAB's \texttt{roots} routine and select the nonzero real root closest to $1$. Since the routine computes the eigenvalues of the companion matrix, the energy correction requires an eigenvalue problem of dimension at most three at each step.

\textbf{Experimental setup.}
We use a Fourier spectral discretization with $N=256$ modes on $N$ equispaced grid points in $[0,2\pi)$. The smooth initial profile, before normalization, is
\begin{equation}\label{eq:initial_smooth}
u_0(x)=\frac{\sin(x)}{2+\cos(x)}.
\end{equation}
For the rough initial data, we prescribe the Fourier coefficients before normalization as
\begin{equation}\label{eq:initial_rough}
\widehat{u}_{0,k}=c_{N,k}\langle k\rangle^{-s},\qquad s>1/2,
\end{equation}
where $\mathbf c_N=(c_{N,k})_k$ is a single draw from the uniform distribution on $[0,1]^N+i[0,1]^N$. We apply the same normalization as for the smooth profile. We use $H^s$ as a label for the data generated with exponent $s$ in \eqref{eq:initial_rough}. %The Fourier truncations are smooth. For an untruncated series in one dimension, this decay with bounded rescaled coefficients guarantees $H^r$ regularity for $r<s-1/2$.

\begin{figure}[htp]
    \centering
    \includegraphics[width=0.75\linewidth]{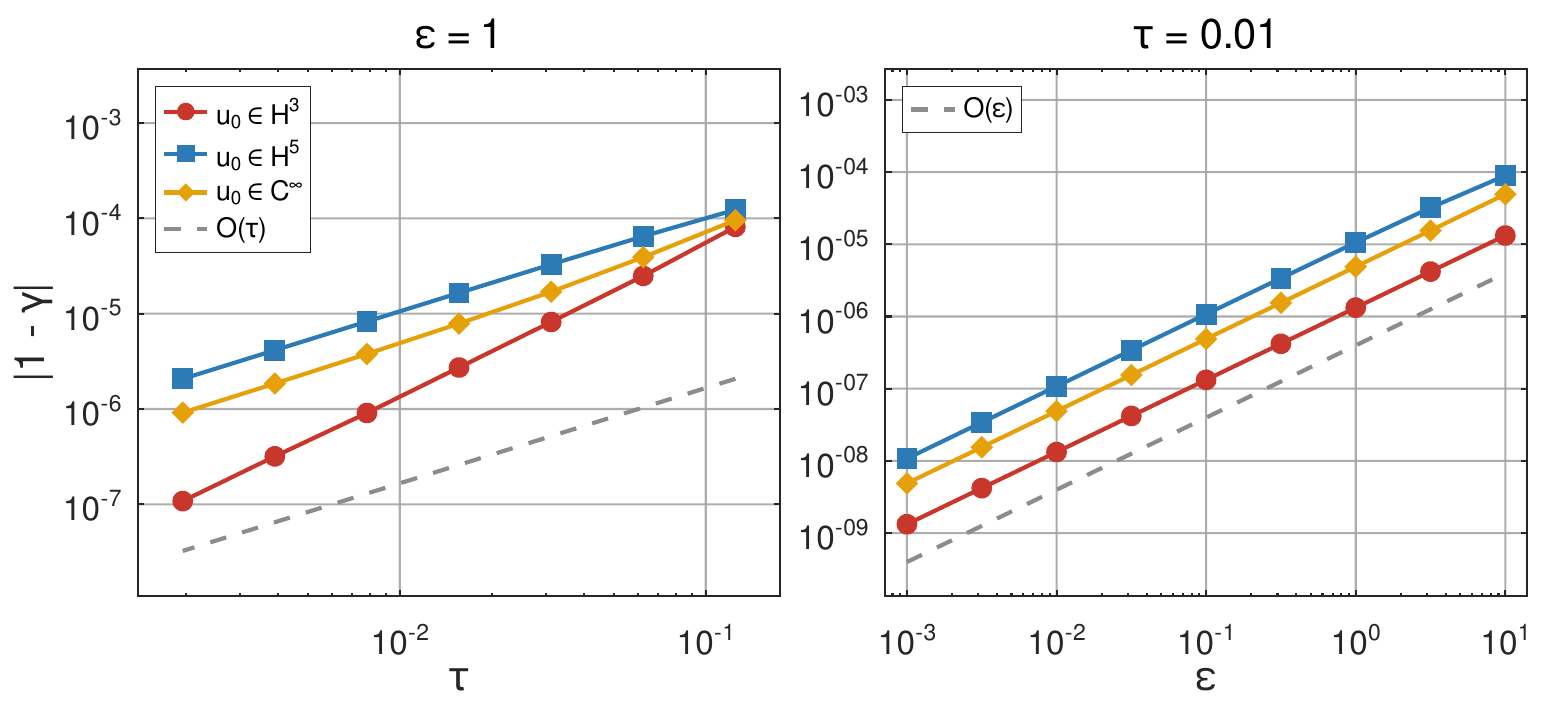}
    \includegraphics[width=0.75\linewidth]{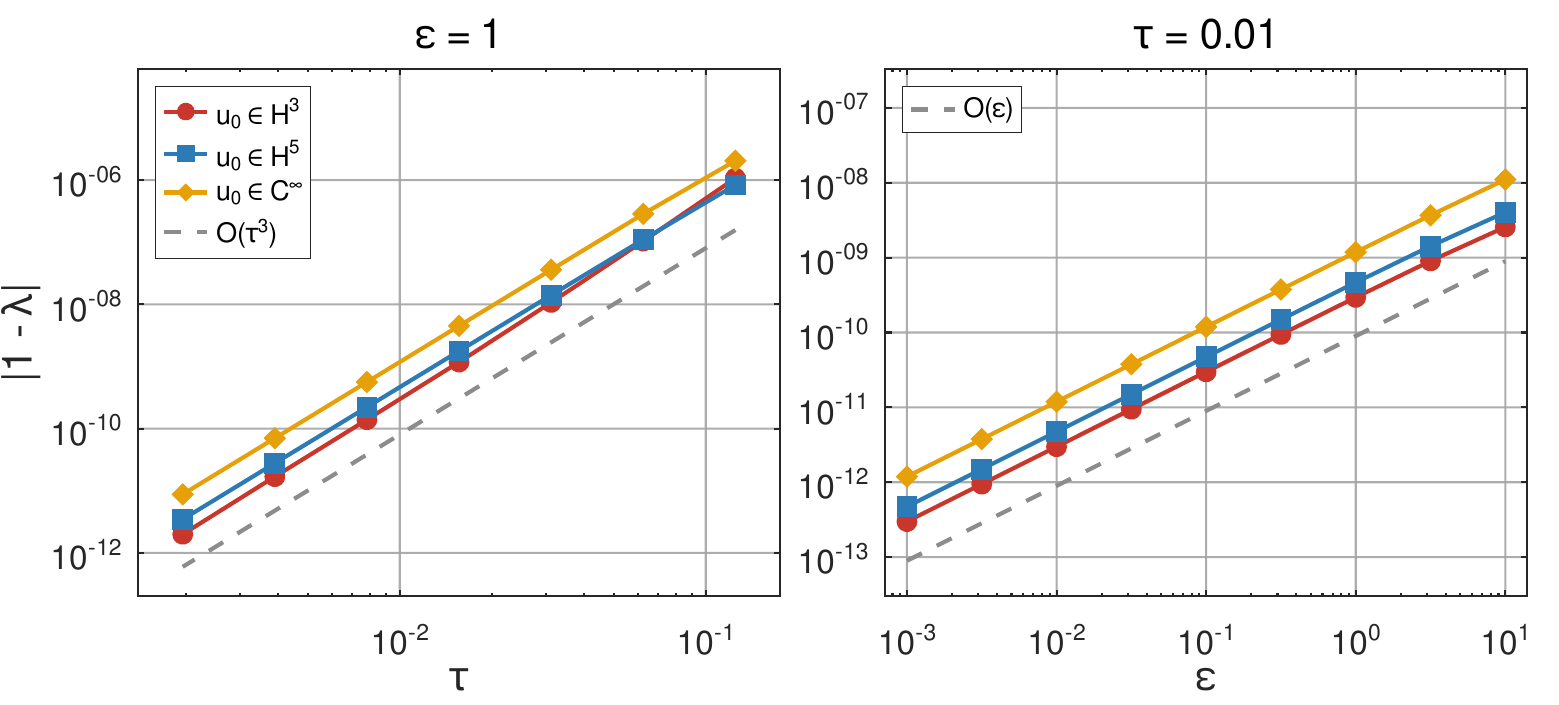}
    \caption{Deviations $|1-\gamma_n|$ (top row) and $|1-\lambda_n|$ (bottom row) for TRLRI2 with the indicated initial data. The left column varies $\tau$ at $\varepsilon=1$. The right column varies $\varepsilon$ at $\tau=0.01$.}
    \label{fig:gamma_lambda}
\end{figure}

The reference solutions are computed with LRI2~\cite{LRI2} using the time step $2\times10^{-5}$. We set $\varepsilon=1$ except when varying $\varepsilon$ in the parameter study in Figure~\ref{fig:gamma_lambda}. The runs at $\varepsilon=1$ test the methods beyond the small-$\varepsilon$ setting of the analysis. All experiments were implemented in MATLAB R2025a and run on a MacBook Pro with a 10-core Apple M4 processor and 24~GB of memory.

\textbf{Convergence of the correction parameters.}
Figure~\ref{fig:gamma_lambda} shows how $|1-\gamma_n|$ and $|1-\lambda_n|$ depend on the time step $\tau$ and the nonlinearity coefficient $\varepsilon$ for one step TRLRI2 starting from a datum $u^n$. For the $H^2$ and smooth data, the step-size study gives rates close to one for $|1-\gamma_n|$ and three for $|1-\lambda_n|$. At fixed $\tau$, both deviations decay linearly with $\varepsilon$. These observations are consistent with the bounds \eqref{eq:gamma2final} and \eqref{eq:lambdaFinal} for $p=2$. For the $H^3$ data, both deviations decay faster with respect to $\tau$ over the tested range.

\begin{remark}
The faster decay of $|1-\gamma|$ in Figure~\ref{fig:gamma_lambda} can be explained by the contribution of high Fourier modes to $\Sigma(\tau)$. For the $H^3$ data, let $m=\|\nabla u_0\|_2^2/\|u_0\|_2^2$ and write
$$
\Sigma(\tau)=\sum_k S_k(\tau),
\qquad
S_k(\tau)=
\bigl(1-\cos(\tau|k|^2)\bigr)(|k|^2-m)|\widehat u_{0,k}|^2.
$$
Consider the modes $B_\tau=\{k\in\mathbb Z:1\le\tau\lvert k\rvert^2\le2\}$. For $\tau m\le 1/2$, these modes satisfy $|k|^2-m\ge |k|^2/2$. Using $\widehat u_{0,k}=a_N c_{N,k}|k|^{-3}$ for $k\ne0$, where $a_N>0$ is the normalization factor, we obtain
$$
S_k(\tau)
\ge \frac{(1-\cos 1)a_N^2}{2}|k|^{-4}|c_{N,k}|^2
\ge \frac{(1-\cos 1)a_N^2}{8}\tau^2|c_{N,k}|^2,
\qquad k\in B_\tau.
$$
If all modes in $B_\tau$ are retained and the average of $|c_{N,k}|^2$ over $B_\tau$ is bounded below, then the count $\#B_\tau\ge c\tau^{-1/2}$ gives
$$
\sum_{k\in B_\tau}S_k(\tau)\ge c_1\tau^{3/2}.
$$
The negative terms satisfy
$$
-\sum_{|k|^2<m}S_k(\tau)
\le \frac{\tau^2}{2}\sum_{|k|^2<m}
|k|^4(m-|k|^2)|\widehat u_{0,k}|^2
=C_0\tau^2.
$$
Thus $\Sigma(\tau)\ge c_1\tau^{3/2}-C_0\tau^2$. When $C_0\sqrt{\tau}\le c_1/2$, the parameter estimate for $p=2$ gives
$$
|1-\gamma|
\lesssim \frac{\varepsilon\tau^3}{\Sigma(\tau)}
\lesssim \varepsilon\tau^{3/2}.
$$
This can account for the faster decay of $|1-\gamma|$ over the tested range. The mass identity \eqref{eq:massidentity} then gives
$$
|1-\lambda|
\lesssim |\delta_M|+\tau^2|1-\gamma|
\lesssim |\delta_M|+\varepsilon\tau^{7/2}.
$$
The faster decay of $|1-\lambda|$ shown in Figure~\ref{fig:gamma_lambda} is then due to the joint effect of $|1-\gamma|$ and $\delta_M$.
\end{remark}

\begin{figure}[htp]
    \centering
    \includegraphics[width=0.95\linewidth]{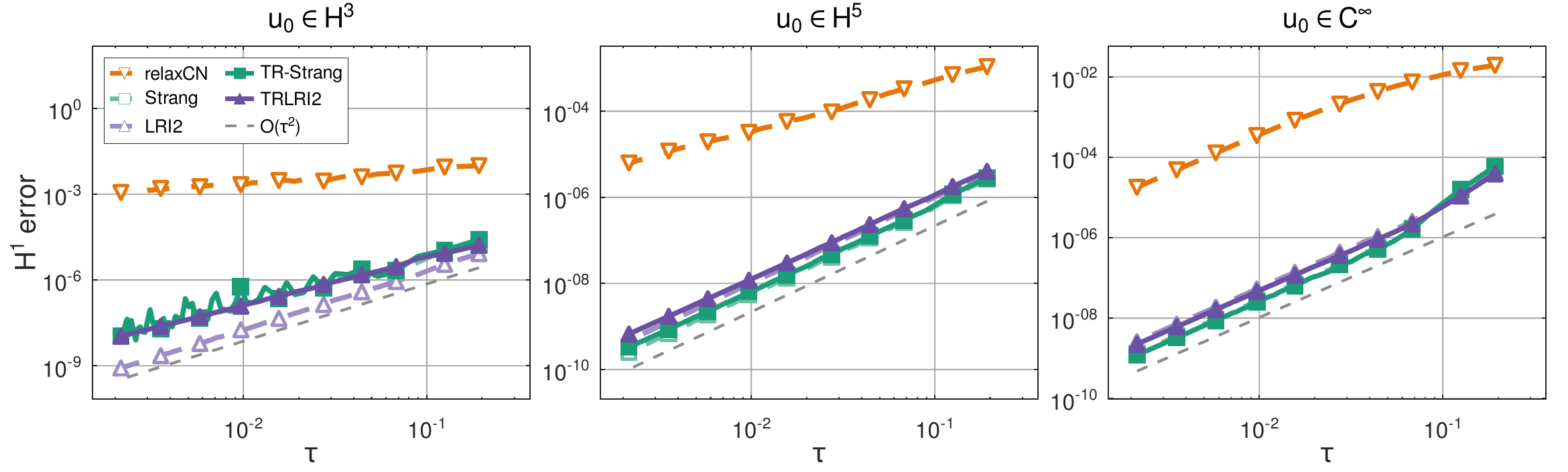}
    \caption{Temporal convergence of the five integrators for the indicated initial data. The $H^1$ error is measured at $T=1$ against the LRI2 reference solution.}
    \label{fig:converge}
\end{figure}

\textbf{Convergence in the $H^1$ norm.}
Figure~\ref{fig:converge} reports the $H^1$ errors of all five methods at $T=1$. For smooth initial data, all methods exhibit second-order convergence. For the $H^5$ data, the observed order of relaxCN falls below two, while the other methods retain second order. For the $H^3$ data, the observed rate of relaxCN decreases further, and the errors of Strang splitting and TRStrang show no clear convergence. LRI2 retains second-order convergence, while the observed order of TRLRI2 is approximately $1.6$.

For $p=2$, inserting the parameter bounds of Theorem~\ref{thm:gamma} into the lower-regularity global error estimate \eqref{eq:globalerror2low} yields a first-order bound under the hypotheses of Theorem~\ref{thm:globalerror2}. The bound in \eqref{eq:globalerror2low} also quantifies how faster decay of the parameter deviations $|1-\gamma_k|$ and $|1-\lambda_k|$ can reduce the accumulated error. The observed order of approximately $1.6$ is consistent with the faster decay seen in Figure~\ref{fig:gamma_lambda} for the $H^3$ data.

\begin{figure}[htp]
    \centering
    \includegraphics[width=0.95\linewidth]{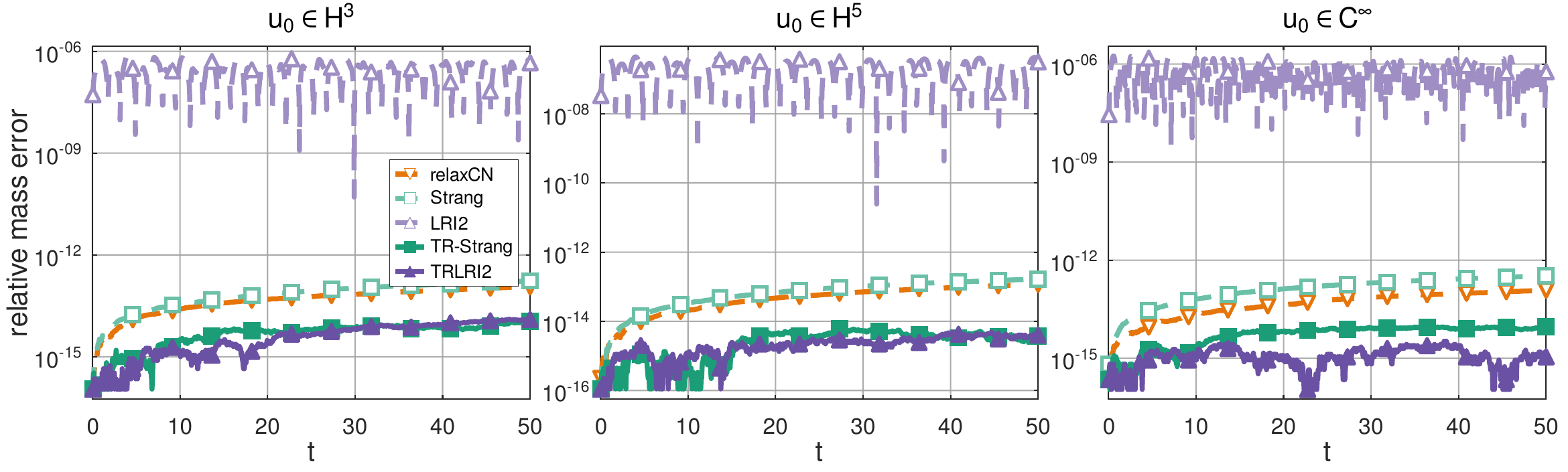}
    \includegraphics[width=0.95\linewidth]{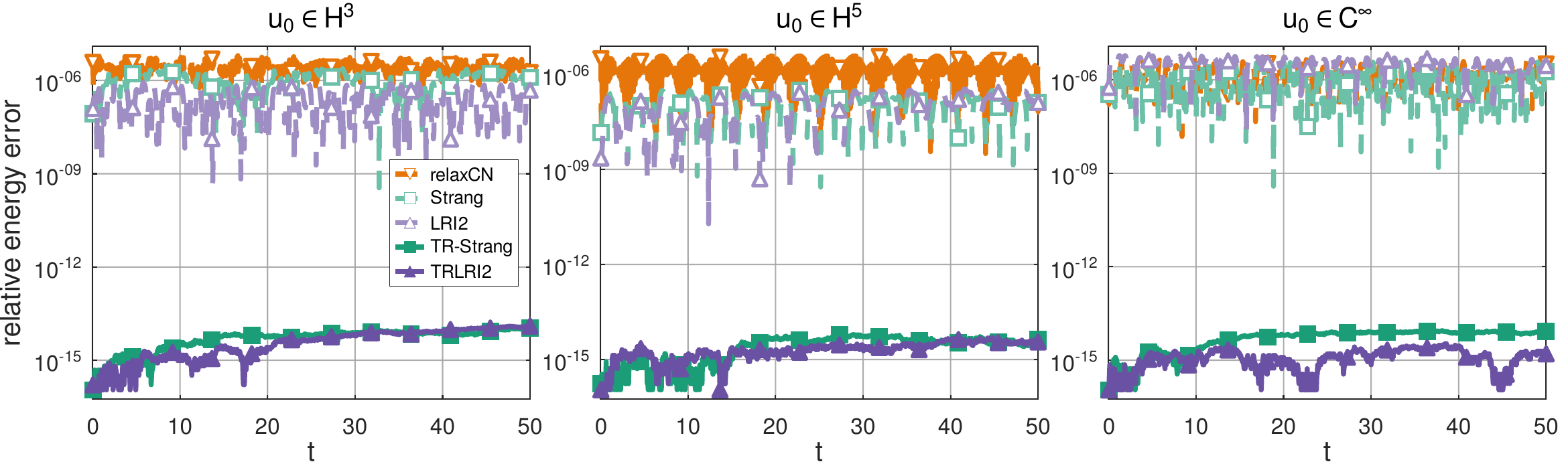}
    \caption{Relative deviations in mass (top row) and physical energy (bottom row) from their initial values for the indicated initial data. The runs extend to $T=50$ with nominal time step $\tau=0.05$.}
    \label{fig:mass-energy}
\end{figure}

\textbf{Conservation of mass and energy.}
We measure the relative deviations in mass and energy by
\[
\frac{\bigl|\mathcal M\left(u^n\right)-\mathcal M\left(u^0\right)\bigr|}{\mathcal M\left(u^0\right)}
\qquad\text{and}\qquad
\frac{|\E(u^n)-\E(u^0)|}{\E(u^0)},
\]
respectively. Figure~\ref{fig:mass-energy} shows these deviations for runs up to $T=50$ with nominal time step $\tau=0.05$. Strang splitting, relaxCN, TRStrang, and TRLRI2 preserve mass in exact arithmetic, whereas LRI2 has no exact mass conservation law. In the experiments, the mass deviations for TRStrang and TRLRI2 are at the level of machine precision and are smaller than those for Strang splitting and relaxCN. The relative physical-energy deviations of Strang splitting, LRI2, and relaxCN are about $10^{-6}$. Both time-relaxed methods preserve the physical energy to machine precision.

\textbf{Computational cost.}
We compare the CPU time required by TRLRI2 and relaxCN to reach a given $H^1$ accuracy. In relaxCN, the system
\[
\left(I-\frac{i\tau}{2}\Delta+\frac{i\tau\varepsilon}{2}\phi^{n+1/2}\right)u^{n+1}
=\left(I+\frac{i\tau}{2}\Delta-\frac{i\tau\varepsilon}{2}\phi^{n+1/2}\right)u^n
\]
is solved by a fixed-point iteration with preconditioner $(I-i\tau\Delta/2)^{-1}$, which is diagonal in Fourier space and avoids assembly and factorization of a dense $N\times N$ matrix. The iteration has contraction factor bounded by $\tau\varepsilon\|\phi^{n+1/2}\|_{L^\infty}/2$ and in practice this factor should be smaller than one. Each iteration requires two FFTs, and we iterate to machine precision. Thus, for a bounded iteration count, it has an $\mathcal O(N\log N)$ cost per step.

\begin{figure}[htp]
    \centering
    \includegraphics[width=0.95\linewidth]{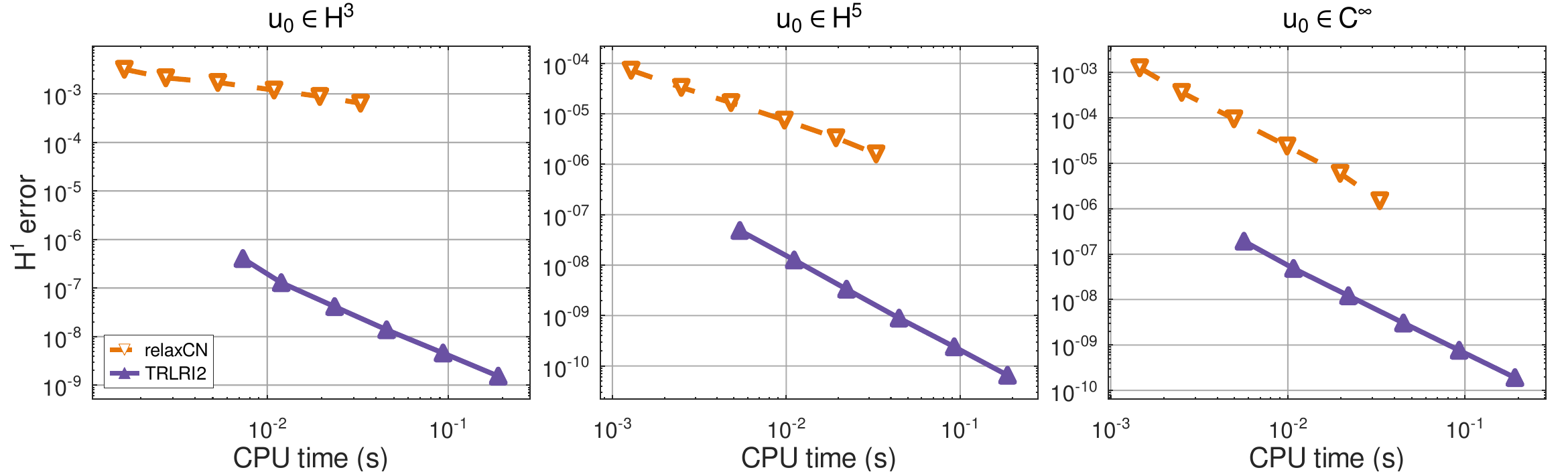}
    \caption{Efficiency comparison of TRLRI2 and relaxCN for the indicated initial data. The $H^1$ error at $T=1$ is plotted against CPU time.}
    \label{fig:efficiency}
\end{figure}

Figure~\ref{fig:efficiency} shows that TRLRI2 (the purple curve with upward-pointing triangle markers) attains a smaller error than relaxCN (the orange dashed curve with downward-pointing triangle marker) at comparable CPU times for all three regularity regimes. In the $H^3$ test, the $H^1$ error of TRLRI2 is approximately $10^{-4}$ times that of relaxCN at comparable CPU times.

\section{Conclusions}\label{sec:conclusions}
In this work, we introduced a novel framework for the construction of mass- and energy-preserving low-regularity integrators for the nonlinear Schrödinger equation. Our approach combines mass normalisation with scalar time relaxation to achieve the simultaneous preservation of both invariants, while retaining favourable convergence properties in the low-regularity regime. We established the existence and uniqueness of the nontrivial relaxation parameter under suitable regularity and spectral nondegeneracy assumptions, together with rigorous convergence estimates for the resulting schemes. Importantly, the proposed correction is computationally inexpensive, requiring only the solution of a polynomial equation of degree at most three. Our theoretical analysis is complemented by numerical experiments demonstrating the exact preservation of mass and energy up to machine precision, alongside favourable accuracy and computational efficiency, including for rough solutions.

A particularly appealing feature of our construction is that it applies to a broad class of admissible baseline low-regularity integrators (cf. Assumption~\ref{assump:lri}), rather than relying on the specific structure of a single numerical scheme. This suggests the possibility of developing more general frameworks for the preservation of multiple invariants in low-regularity numerical integration. In particular, the framework extends to the Korteweg–de Vries equation and the investigation of the compatibility of this construction with general frameworks for resonance-based integration \cite{LRI2,maierhofer2023bridging}, provides an interesting direction for future research.
\section*{Acknowledgments} 
The work of TJ is supported by HKUST Postgraduate Studentship, Overseas Research Award (HKUST, March - August, 2024), and the European Union (ERC, COCOA, 101170147). 
The authors used ChatGPT (models 6-astra and 6.1-sol) to assist with implementing the comparison method relaxCN, generate the TikZ code for Figure~\ref{fig:comparison_diagram}, and proofread the manuscript, and adopted its suggestion to use the auxiliary set $J_\tau$ in the proof of Theorem~\ref{thm:global}.

\bibliographystyle{plain} % alphabetical order
\bibliography{references}
\end{document}